\documentclass[11pt,a4paper]{article}

\usepackage{amssymb}
\usepackage{dsfont}
\usepackage{amsmath}
\usepackage{amsthm}
\usepackage[dvipsnames]{xcolor}
\usepackage{hyperref}
\hypersetup{
    colorlinks=true,
    linkcolor=Green,
    citecolor=Green,
    urlcolor=Cyan
    }
    
\usepackage{soul}
\usepackage{diagbox}
\usepackage{pgfplotstable}
\usepackage{pgfplots}
\usepgfplotslibrary{fillbetween}
\pgfplotsset{compat = newest}
\usepackage{multirow}
\pgfplotsset{compat=1.18}
\usepackage{subcaption}
\usepackage[titletoc,title]{appendix}
\usepackage[ruled]{algorithm2e}
\usepackage{nomencl}
\usepackage{etoolbox}
\renewcommand\nomgroup[1]{%
  \item[\bfseries
  \ifstrequal{#1}{A}{Abbreviations}{%
  \ifstrequal{#1}{B}{Roman symbols}{%
  \ifstrequal{#1}{C}{Greek symbols}{%
  \ifstrequal{#1}{D}{Subscripts}{%
  \ifstrequal{#1}{E}{Superscripts}{%
  \ifstrequal{#1}{F}{Other symbols}{}}}}}}%
]}
\usepackage{tikz,float}
\usetikzlibrary{positioning, fit, calc}   
\tikzset{block/.style={draw, thick, text width=1.25cm, minimum height=0.8cm, align=center},   
line/.style={-latex}}
\pgfplotsset{every axis/.append style={
    label style={font=\scriptsize},
    tick label style={font=\scriptsize}  
    }}

\usepackage{titlesec}

\titleformat{\paragraph}
{\normalfont\normalsize\bfseries}{\theparagraph}{1em}{}
\titlespacing*{\paragraph}
{0pt}{3.25ex plus 1ex minus .2ex}{1.5ex plus .2ex}

\def\num#1{\numx#1}\def\numx#1e#2{{#1}\mathrm{e}{#2}}
\newcommand*{\collauthor}[2]{{#1}$^{#2}$}

\newcommand*{\affiliation}[2]{$\mbox{}^{{#2}}${#1}}
\newcommand*{\colltitle}[1]{\textbf{#1}}

\newenvironment{keywords}[1]{\vspace{1cm}\\{\bf \slshape{Keywords}}\quad\slshape{#1}}{}
\DeclareMathOperator*{\argmin}{arg\,min}
\newtheorem{exmp}{Example}[]

\DeclareMathOperator{\Tr}{Tr}
\DeclareMathOperator{\Hessian}{Hess}
\newcommand*\Bell{\ensuremath{\boldsymbol\ell}}

\newcommand{\ifrac}{\mathfrak{i}}
\newcommand{\jfrac}{\mathfrak{j}}
\newcommand{\Kf}{\mathfrak{K}}
\newcommand{\nf}{\mathfrak{n}}

\begin{document}
\begin{center}
\begin{Large}
  \colltitle{Uncertainty quantification for deep learning-based schemes for solving high-dimensional backward stochastic differential equations}
\end{Large}
\vspace*{1.5ex}

\begin{sc}
\begin{large}
\collauthor{Lorenc Kapllani}{}, \collauthor{Long Teng}{} and \collauthor{Matthias Rottmann}{}
\end{large}
\end{sc}
\vspace{1.5ex}

\affiliation{Chair of Applied and Computational Mathematics,\\Faculty of Mathematics and Natural Sciences,\\
University of Wuppertal,\\
Gau{\ss}str. 20, 42119 Wuppertal, Germany\linebreak }{} \\
\end{center}

\section*{Abstract}
Deep learning-based numerical schemes for solving high-dimensional backward stochastic differential equations (BSDEs) have recently raised plenty of scientific interest. While they enable numerical methods to approximate very high-dimensional BSDEs, their reliability has not been studied and is thus not understood. In this work, we study uncertainty quantification (UQ) for a class of deep learning-based BSDE schemes. More precisely, we review the sources of uncertainty involved in the schemes and numerically study the impact of different sources. Usually, the standard deviation (STD) of the approximate solutions obtained from multiple runs of the algorithm with different datasets is calculated to address the uncertainty. This approach is computationally quite expensive, especially for high-dimensional problems. Hence, we develop a UQ model that efficiently estimates the STD of the approximate solution using only a single run of the algorithm. The model also estimates the mean of the approximate solution, which can be leveraged to initialize the algorithm and improve the optimization process. Our numerical experiments show that the UQ model produces reliable estimates of
the mean and STD of the approximate solution for the considered class of deep learning-based BSDE schemes. The estimated STD captures multiple sources of uncertainty, demonstrating its effectiveness in quantifying the uncertainty. Additionally, the model illustrates the improved performance when comparing different schemes based on the estimated STD values. Furthermore, it can identify hyperparameter values for which the scheme achieves good approximations.
\begin{keywords}
backward stochastic differential equations, high-dimensional problems, deep neural networks, uncertainty quantification, heteroscedastic nonlinear regression
\end{keywords}



\section{Introduction}
\label{sec1}
Deep learning has attracted the interest of academics in various fields, including computer vision and natural language processing, as well as traditional sciences such as physics, chemistry, biology, and finance. Since the predictions of models that use deep learning in decision-making processes are prone to noise and model errors~\cite{jiang2018trust}, assessing the model's reliability before it can be used in practice is critical. An example of such decisions is the pricing and hedging of different contracts in finance. Companies may suffer from significant financial losses as a result of poor judgments. Thus, it is highly desirable to understand the uncertainties in deep learning-based numerical schemes and develop methods to quantify them.

Backward stochastic differential equations (BSDEs) are important tools used to model problems in scientific fields due to their connections to partial differential equations (PDEs) and stochastic control problems via the nonlinear Feynman–Kac formula~\cite{El1997}. In most cases, nonlinear BSDEs cannot be solved explicitly. Hence, advanced numerical techniques to approximate their solutions become desired. Over the years, many numerical methods have been proposed for solving BSDEs, e.g., \cite{bouchard2004discrete,zhang2004numerical,gobet2005regression,lemor2006rate,zhao2006new,bender2008,ma2008numerical,zhao2010stable,gobet2010solving,crisan2012solving,zhao2014new,ruijter2015fourier,ruijter2016fourier}. However, most of those methods are not suitable for solving high-dimensional BSDEs due to the curse of dimensionality. Some techniques such as parallel computing using GPU computing~\cite{gobet2016stratified, kapllani2022multistep} or sparse grid methods~\cite{zhang2013sparse,fu2017efficient,chassagneux2022learning} can only solve moderate dimensional BSDEs within a reasonable computational time. 

Recently, the authors in~\cite{weinan2017deep,han2018solving} developed a deep learning-based scheme (we refer to it as the DBSDE scheme in the rest of the paper) for solving high-dimensional nonlinear BSDEs. The method employs deep neural networks (DNNs) to approximate the unknown gradient of the solution by reformulating the BSDE problem as a stochastic optimization problem. Due to the universal approximation capability~\cite{hornik1989multilayer,cybenko1989approximation} of neural networks (NNs), the objective function can be effectively optimized in practice. Therefore, the function values of interest (the unknown solution and its gradient) are obtained quite accurately. The method has opened the door to solve such problems in hundreds of dimensions within a reasonable computational time. After the DBSDE scheme, there are several works proposed to improve it, e.g., by algorithmic adjustments or methodological extensions~\cite{raissi2018forward,fujii2019asymptotic,hure2020deep,beck2021deep,takahashi2022new,germain2022approximation,kapllani2022effect,kapllani2020deep}. Furthermore, some convergence analysis of the DBSDE scheme have also been studied, e.g., see~\cite{han2020convergence,jiang2021convergence} for the error analysis (utilizing universal approximation capability of NNs) and  \cite{wang2021gradient} for its gradient convergence (under a restrictive choice of NN setting).

As an example, we consider the pioneering DBSDE scheme to study uncertainty quantification (UQ) and introduce our UQ model. In the DBSDE scheme, the authors reformulate the BSDE as a stochastic control problem and use the Euler-Maruyama method to discretize the integrals. The unknown functions (solution at initial time and its gradient on the whole time domain) are estimated using DNNs. The parameters of DNNs are then optimized using the stochastic gradient descent (SGD) algorithm. The DBSDE method incorporates various sources of uncertainty, including finite time discretization, restrictive choice of DNN specifications, the lack of convergence guarantees of the SGD algorithm, and finite sample size during stochastic optimization. It is crucial to identify and quantify these different sources of uncertainty in the DBSDE scheme for practical applications. Therefore, we review these sources of uncertainty and numerically demonstrate the impact of different sources. Our numerical experiments show that it is practically challenging to disentangle them, emphasizing the importance of quantifying uncertainty before using the scheme in practice. Usually, the standard deviation (STD) of the approximate solutions obtained from multiple runs of the DBSDE algorithm with different datasets is calculated to account for the uncertainty in a given prediction, see~\cite{weinan2017deep}. This approach is computationally expensive, especially in high-dimensional cases. Therefore, we develop a UQ model to estimate the STD of the approximate solution without requiring multiple runs. Several techniques have been proposed in the literature to quantify uncertainty, such as Monte Carlo Dropout~\cite{gal2016dropout}, Monte Carlo DropConnnect~\cite{mobiny2021dropconnect}, deep
ensembles~\cite{lakshminarayanan2017simple,ovadia2019can}, Flipout-based variational inference~\cite{wen2018flipout}, Markov Chain Monte Carlo~\cite{kupinski2003ideal}, and many others~\cite{hendrycks2016baseline,oberdiek2018classification,oberdiek2022uqgan}. A review of these techniques can be found in~\cite{abdar2021review}. To the best of our knowledge, there are no applications or developments of UQ models specifically for deep learning-based BSDE schemes. Hence, we develop a UQ model with the aim of addressing this gap.

Our UQ model is based on a commonly used approach for quantifying uncertainty in heteroscedastic nonlinear regression, see~\cite{beal1988heteroscedastic} for heteroscedastic least square type regression methods and~\cite{nix1994estimating,lakshminarayanan2017simple,kendall2017uncertainties} for heteroscedastic NN regression methods. We make the assumption that the residuals or errors of the DBSDE scheme follow a normal distribution with zero mean and the STD depending on the chosen parameter set of the discretized BSDE. This is a standard assumption in heteroscedastic regression. In our method, we use a DNN to learn two functions that estimate the mean and STD of the approximate solution. To train the DNN, we construct a dataset of independent and identically distributed (i.i.d) samples, which includes different parameter sets of the discretized BSDE and the corresponding approximate solutions obtained from the DBSDE algorithm. After generating a moderate number of samples, we optimize the network parameters by minimizing the negative log-likelihood. Our UQ model provides an estimate of the STD of the approximate solution in a more computationally efficient manner compared to running multiple iterations of the algorithm per parameter set. Additionally, the estimated mean of the approximate solution from our model can be used to initialize the algorithm and improve the optimization process. Note that the proposed UQ model is applicable not only to the DBSDE scheme but also to other deep learning-based schemes for solving BSDEs. In addition to the DBSDE scheme, we apply our UQ model to the Locally additive DBSDE (LaDBSDE) scheme~\cite{kapllani2020deep}, which exhibits better convergence than the DBSDE scheme. Our numerical experiments demonstrate that the UQ model produces reliable estimates of the mean and STD of the approximate solution for both the DBSDE and LaDBSDE schemes. Moreover, we show multiple practical implications of using the UQ model. Firstly, the estimated STD captures multiple sources of uncertainty, demonstrating its effectiveness in quantifying the uncertainty. Secondly, the UQ model illustrates the improved performance of the LaDBSDE scheme in comparison to the DBSDE scheme based on the corresponding estimated STD values. Finally, our UQ model can be utilized to determine DNN hyperparameter values for which the respective scheme performs well, e.g.\ the number of discretization points.

The remainder of the paper is organized as follows. In Section~\ref{sec2}, we provide the necessary foundations to view the BSDE as a learning problem. In Section~\ref{sec3}, we introduce the DBSDE scheme for solving high-dimensional BSDEs, discuss the sources of uncertainty in the scheme, and present a UQ model to estimate the STD of the approximate solution. In Section~\ref{sec4}, we analyze the practical impact of different sources of uncertainty in the DBSDE scheme and demonstrate the effectiveness of our UQ model through numerical tests for both the DBSDE and LaDBSDE schemes. Finally, in Section~\ref{sec5}, we conclude our work.

\section{Backward stochastic differential equations as a learning problem}
\label{sec2}
In the following, we introduce the notions of the BSDEs that are used throughout the paper.

\subsection{Preliminaries}
\label{subsec21}
Let $\left(\Omega,\mathcal{F},\mathbb{P},\{\mathcal{F}_t\}_{0\le t \le T}\right)$ be a complete, filtered probability space. In this space a standard $d$-dimensional Brownian motion $W_t$ is defined, such that the filtration $\{\mathcal{F}_t\}_{0\le t\le T}$ is the natural filtration of $W_t.$ Throughout the whole paper we rely on the following notations
\begin{itemize}
    \item $| x |$ for the standard Euclidean norm of $x \in \mathbb{R}$ or $x \in \mathbb{R}^{d}$.
    \item $\mathbb{S}^2\left( [0, T] \times \Omega; \mathbb{R}^{d} \right)$ for the space of continuous and progressively measurable stochastic processes $X: [0, T] \times \Omega \to \mathbb{R}^{d} $ such that $\mathbb{E}\bigl[\sup_{ 0 \leq t\leq T}|X_t|^2\bigr] < \infty$.
    \item  $\Delta = \{t_n|t_n \in [0, T], n = 0, 1, \ldots, N, t_n < t_{n+1}, \Delta t = t_{n+1} - t_{n}, t_0 = 0, t_N = T\}$ is a uniform discretization of the time interval $[0, T]$.
    \item $\mathbb{H}^{2}\left( [0, T] \times \Omega ; \mathbb{R}^{1 \times d} \right)$ for the space of progressively measurable stochastic processes $Z: [0, T] \times \Omega \to \mathbb{R}^{1 \times d} $ such that $\mathbb{E}\bigl[\int_0^T |Z_t|^2 \, dt\bigr] < \infty$.    
    \item $\mathbb{H}^{\Delta, 2}\left( \{0, 1, \ldots, N\} \times \Omega ; \mathbb{R}^{1 \times d} \right)$ for the space of progressively measurable discrete stochastic processes $Z^{\Delta}: \{0, 1, \ldots, N\} \times \Omega \to \mathbb{R}^{1 \times d} $ such that $\mathbb{E}\bigl[\sum_{n=0}^{N} |Z_{t_n}^{\Delta}|^2 \, \Delta t\bigr] < \infty$.
    \item $\mathbb{L}^2_{\mathcal{F}_t}\left(\Omega ; \mathbb{R}^{d} \right)$ for the space of $\mathcal{F}_t$-measurable random variables $\xi: \Omega \to \mathbb{R}^{d} $ such that $\mathbb{E}\bigl[|\xi|^2 \bigr] < \infty$.
    \item $C^{1,2}\left( [0, T] \times \mathbb{R}^{d}; \mathbb{R} \right)$ for the space of continuous and differentiable functions of two arguments $u: [0, T] \times \mathbb{R}^{d} \to \mathbb{R}$, differentiable with respect to the first argument and twice differentiable with respect to the second argument.
    \item $D^{\top}$ for the transpose of matrix $D \in \mathbb{R}^{d \times d}$.
    \item $\Tr\left[D\right]$ for the trace of matrix $D \in \mathbb{R}^{d \times d}$.
    
\end{itemize}
\subsection{Nonlinear Feynman-Kac formula}
\label{subsec22}
We present the nonlinear Feynman–Kac formula, which provides a probabilistic representation for the solution of a semi-linear parabolic PDE through a (decoupled) forward-backward stochastic differential equation (FBSDE).

Let $T \in (0, \infty)$, $d \in \mathbb{N}$, the functions $f: [0, T ] \times \mathbb{R}^d \times \mathbb{R} \times \mathbb{R}^{1 \times d} \to \mathbb{R}$ and $g: \mathbb{R}^d \to \mathbb{R}$ are continuous, and $(a, b): [0, T ] \times \mathbb{R}^d \to \mathbb{R}^d \times \mathbb{R}^{d\times d}$ is continuously differentiable for all variables. Let $u \in C^{1,2}([0, T] \times \mathbb{R}^d; \mathbb{R})$ satisfy that $u(T, x) = g(x)$ and the semi-linear parabolic PDE
\begin{equation}
    \frac{\partial u}{\partial t} + \nabla u(t, x)\,a(t, x) + \frac{1}{2} \Tr\left[b b^{\top} \Hessian u (t, x)\right] + f\left(t, x, u, \nabla u \,  b\right)(t, x) = 0,
\label{eq1}
\end{equation}
for all $(t, x) \in ([0, T]\times \mathbb{R}^d)$, where $\Hessian u$ is the Hessian matrix and $\nabla u$ the gradient of $u$ with respect to the spatial variable $x$. Consider the following decoupled FBSDE
\begin{equation}
    \begin{split}
        \left\{
            \begin{array}{rcl}
                X_t & = & x_0 + \int_{0}^{t} a \left(s, X_s\right)\,ds + \int_{0}^{t} b \left(s, X_s\right)\,dW_s,\\
   	   		    Y_t & = & g\left(X_T\right) + \int_{t}^{T} f\left(s, X_s, Y_s, Z_s\right)\,ds -\int_{t}^{T} Z_s \,dW_s,
            \end{array}
        \right. 
    \end{split}
\label{eq2}
\end{equation}
where $f:\left[0, T\right]\times \mathbb{R}^d\times \mathbb{R}\times\mathbb{R}^{1\times d} \to \mathbb{R}$ is the driver function, $f(t, X_t, Y_t, Z_t)$ is $\mathcal{F}_t$-adapted, the third term on the right-hand side is an It\^o-type integral and $g:\mathbb{R}^d \to \mathbb{R}$ is the function for the terminal condition depending on the final value $X_T$ of the forward stochastic differential equation (SDE). The triple of processes $\{\left(X_t,Y_t,Z_t\right)\}_{0\leq t \leq T} \in \mathbb{S}^2\left( [0, T] \times \Omega; \mathbb{R}^d \right) \times \mathbb{S}^2\left( [0, T] \times \Omega; \mathbb{R} \right) \times \mathbb{H}^2\left( [0, T] \times \Omega; \mathbb{R}^{1 \times d} \right)$ is the solution of the above FBSDE if it satisfies~\eqref{eq2} $\mathbb{P}$-a.s. $\forall$ $t \in \left[0,T\right]$. Assume that \eqref{eq1} has a classical solution $u \in C^{1,2}([0, T] \times \mathbb{R}^d; \mathbb{R})$ and the usual regularity conditions of~\eqref{eq2} are satisfied~\cite{El1997}. 
Then the solution of~\eqref{eq2} can be represented $\mathbb{P}$-a.s. by 
\begin{equation}
	Y_t = u\left(t, X_t\right), \quad Z_t= \nabla u\left(t, X_t\right)b \left(t, X_t\right) \quad \forall \,\, t \in \left[0,T\right).
\label{eq3}
\end{equation}

\subsection{Formulation of the BSDE as a suitable stochastic control problem}
\label{subsec23}
We now view the solution $(X_t, Y_t, Z_t)$ of~\eqref{eq2} as the stochastic control problem
\begin{align}
    \begin{split}
        &\inf_{Y_0 \in \mathbb{L}^2_{\mathcal{F}_0}\left(\Omega; \mathbb{R} \right), \, Z \in \mathbb{H}^2\left( [0, T] \times \Omega;\mathbb{R}^{1 \times d} \right)} \mathbf{L}\left(Y_0, Z\right),\\
         s.t. \quad & X_t = x_0 + \int_{0}^{t} a \left(s, X_s\right)\,ds + \int_{0}^{t} b \left(s, X_s\right)\,dW_s,\\
         & Y_t =  Y_0 - \int_{0}^{t} f\left(s, X_s, Y_s, Z_s\right)\,ds +\int_{0}^{t} Z_s\,dW_s,
    \end{split}
\label{eq4}
\end{align}
for $t \in [0, T]$, where 
$$\mathbf{L}\left(Y_0, Z\right):= \mathbb{E}\bigl[|g(X_T)-Y_T|^2\bigr].$$ 
The solution of~\eqref{eq2} is a minimizer of~\eqref{eq4} since the loss function attains zero when it is evaluated at the solution.  In addition, the wellposedness of the BSDEs (under the usual regularity conditions~\cite{El1997}) ensures the existence and uniqueness of the minimizer. Due to~\eqref{eq3}, we seek a function approximator for $u: \mathbb{R}^d \to \mathbb{R}$ and $\nabla u\,b: [0, T] \times \mathbb{R}^d \to \mathbb{R}^{1 \times d}$ to approximate the unknown solution $Y_0 = u(t_0, x_0)$ and  $Z_t = \nabla u\left(t, X_t\right) b\left(t, X_t\right)$ $\forall$ $t \in [0, T]$. Due to their approximation capability in high dimensions, NNs are a promising candidate.

\subsection{Neural networks as function approximators}
\label{subsec24}
For our purpose, we consider fully connected feedforward NNs or DNNs. Let $d_0, d_1\in \mathbb{N}$ be the input and output dimensions, respectively. We fix the global number of layers as $L+2$, $L \in \mathbb{N}$ the number of hidden layers each with $\eta \in \mathbb{N}$ neurons. The first layer is the input layer with $d_0$ neurons and the last layer is the output layer with $d_1$ neurons. A DNN is a function $\phi(\cdot; \theta): \mathbb{R}^{d_0} \to \mathbb{R}^{d_1}$ composed of a sequence of simple functions, which therefore can be collected in the following form
\begin{equation*}
    x \in \mathbb{R}^{d_0} \longmapsto A_{L+1}(\cdot;\theta(L+1)) \circ \varrho \circ A_{L}(\cdot;\theta(L)) \circ \varrho \circ \ldots \circ \varrho \circ A_1(x;\theta(1)) \in \mathbb{R}^{d_1},
\end{equation*}
where $\theta:=\left( \theta(1), \ldots, \theta(L+1) \right) \in \mathbb{R}^{P}$ and $P$ is the total number of network parameters, $x \in \mathbb{R}^{d_0}$ is called an input vector. Moreover, $A_l(\cdot; \theta(l)), l = 1, 2, \ldots, L+1$ are affine transformations: $A_1(\cdot;\theta(1)): \mathbb{R}^{d_0} \to \mathbb{R}^{\eta}$, $A_l(\cdot;\theta(l)),  l = 2, \ldots, L: \mathbb{R}^{\eta} \to \mathbb{R}^{\eta}$ and $A_{L+1}(\cdot;\theta(L+1)): \mathbb{R}^{\eta} \to \mathbb{R}^{d_1}$, represented by
\begin{equation*}
    A_l(v;\theta(l)) = \mathcal{W}_l v + \mathcal{B}_l, \quad v \in \mathbb{R}^{\eta_{l-1}},
\end{equation*}
where $\theta(l):=\left(\mathcal{W}_l, \mathcal{B}_l\right)$, $\mathcal{W}_l \in \mathbb{R}^{\eta_{l} \times \eta_{l-1}}$ is the weight matrix and $\mathcal{B}_l \in \mathbb{R}^{\eta_{l}}$ is the bias vector with $\eta_0 = d_0, \eta_{L+1} = d_1, \eta_l = \eta$ for $l = 1, \ldots, L$ and $\varrho: \mathbb{R} \to \mathbb{R}$ is a nonlinear function (called the activation function), and applied component-wise on the outputs of $A_l(\cdot;\theta(l))$. Common choices are $\tanh(\cdot), \sin(\cdot), \max(0,\cdot)$ etc. All these matrices $\mathcal{W}_l$ and vectors $\mathcal{B}_l$ form the parameters $\theta$ of the DNN and can be collected as 
$$P = \sum_{l=1}^{L+1}\eta_{l}(\eta_{l-1}+1) = \eta(d_0+1) + \eta(\eta+1)(L-1) + d_1(\eta+1),$$
for fixed $d_0, d_1, L$ and $\eta$. We denote by $\Theta$ the set of possible parameters for the DNN $\phi(\cdot; \theta)$ with $\theta \in \Theta$. The Universal Approximation Theorem~\cite{hornik1989multilayer,cybenko1989approximation} justifies the use of NNs as function approximators.

\section{Uncertainty in the DBSDE scheme}
\label{sec3}
In this section, we discuss the sources of uncertainty in the DBSDE scheme and propose a UQ model to estimate the STD of the approximate solution.

\subsection{The DBSDE scheme}
\label{subsec31}
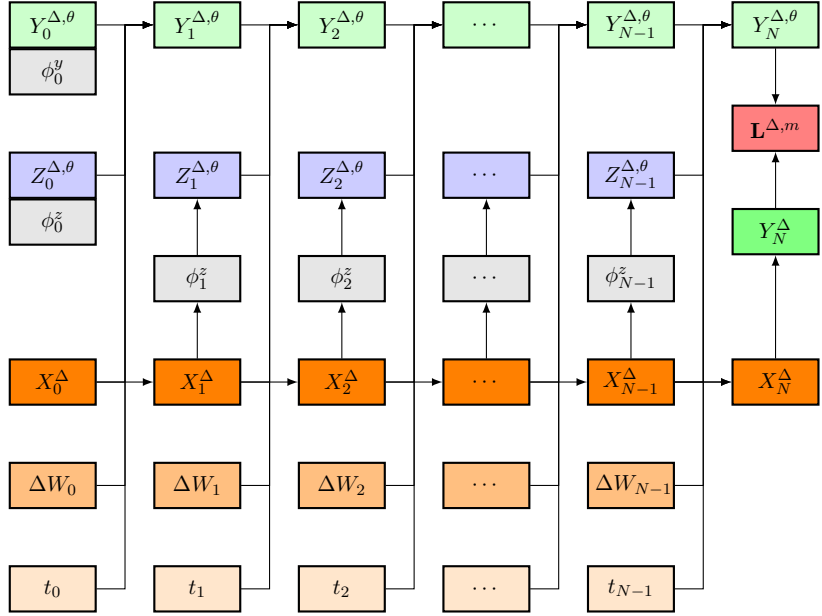
\begin{figure}[htb!]
    \centering
    \begin{tikzpicture}[scale=0.75, transform shape]  
        \node[block, fill=orange!20] (t0) {$t_0$};  
        \node[block,right=of t0, fill=orange!20] (t1) {$t_1$};   
        \node[block,right=of t1, fill=orange!20] (t2) {$t_2$};  
        \node[block,right=of t2, fill=orange!20] (ti) {$\cdots$};
        \node[block,right=of ti, fill=orange!20] (tNm1) {$t_{N-1}$};
        
        \node[block, above = of t0, fill=orange!50] (DW0) {$\Delta W_0$};   
        \node[block,right=of DW0, fill=orange!50] (DW1) {$\Delta W_1$};  
        \node[block,right=of DW1, fill=orange!50] (DW2) {$\Delta W_2$};  
        \node[block,right=of DW2, fill=orange!50] (DWi) {$\cdots$};  
        \node[block,right=of DWi, fill=orange!50] (DWNm1) {$\Delta W_{N-1}$};  
        
        \node[block, above = of DW0, fill=orange] (XD0) {$X^{\Delta}_0$};   
        \node[block,right=of XD0, fill=orange] (XD1) {$X^{\Delta}_1$};  
        \node[block,right=of XD1, fill=orange] (XD2) {$X^{\Delta}_2$};
        \node[block,right=of XD2, fill=orange] (XDi) {$\cdots$};  
        \node[block,right=of XDi, fill=orange] (XDNm1) {$X^{\Delta}_{N-1}$};  
        \node[block,right=of XDNm1, fill=orange] (XDN) {$X^{\Delta}_N$};  
        
        \node[block,above = of XD1, fill=gray!20] (phiz1) {$\phi^z_1$};  
        \node[block,right=of phiz1, fill=gray!20] (phiz2) {$\phi^z_2$};
        \node[block,right=of phiz2, fill=gray!20] (phizi) {$\cdots$};  
        \node[block,right=of phizi, fill=gray!20] (phizNm1) {$\phi^z_{N-1}$};  
        
        \node[block,above = of phiz1, fill=blue!20] (ZDT1) {$Z_1^{\Delta, \theta}$};  
        \node[block,left=of ZDT1, fill=blue!20] (ZDT0) {$Z_0^{\Delta, \theta}$};
        \node[block,below= 0mm of ZDT0, fill=gray!20] (phiz0) {$\phi_0^z$};
        \node[block,right=of ZDT1, fill=blue!20] (ZDT2) {$Z_2^{\Delta, \theta}$};  
        \node[block,right=of ZDT2, fill=blue!20] (ZDTi) {$\cdots$};  
        \node[block,right=of ZDTi, fill=blue!20] (ZDTNm1) {$Z_{N-1}^{\Delta, \theta}$};  
        
        \node[block,above = of ZDT0, fill=gray!20] (phiy0) {$\phi_0^y$};  
        \node[block,above = 0mm of phiy0, fill=green!20] (YDT0) {$Y_0^{\Delta, \theta}$};  
        \node[block,right =of YDT0, fill=green!20] (YDT1) {$Y_1^{\Delta, \theta}$};  
        \node[block,right =of YDT1, fill=green!20] (YDT2) {$Y_2^{\Delta, \theta}$};  
        \node[block,right =of YDT2, fill=green!20] (YDTi) {$\cdots$};  
        \node[block,right =of YDTi, fill=green!20] (YDTNm1) {$Y_{N-1}^{\Delta, \theta}$};  
        \node[block,right = of YDTNm1, fill=green!20] (YDTN) {$Y_N^{\Delta, \theta}$};  
        \node[block,below=of YDTN, fill=red!50] (LDm) {$\mathbf{L}^{\Delta,m}$};  
        \node[block,below=of LDm, fill=green!50] (YDN) {$Y^{\Delta}_N$};

        \draw[line] (XD0)--(XD1);
        \draw[line] (XD1)--(XD2);
        \draw[line] (XD2)--(XDi);
        \draw[line] (XDi)--(XDNm1);
        \draw[line] (XDNm1)--(XDN);
        \draw[line] (XDN)--(YDN);

        \draw[line] (XDNm1)--(XDN);
        \draw[line] (XD1)--(phiz1);
        \draw[line] (XD2)--(phiz2);
        \draw[line] (XDi)--(phizi);
        \draw[line] (XDNm1)--(phizNm1);
        
        \draw[line] (phiz1)--(ZDT1);
        \draw[line] (phiz2)--(ZDT2);
        \draw[line] (phizi)--(ZDTi);
        \draw[line] (phizNm1)--(ZDTNm1);
        
        \draw[line] (YDT0)--(YDT1);
        \draw[line] (t0) -| ($(t0) !0.5! (YDT1)$) |- (YDT1);
        \draw[line] (DW0) -| ($(DW0) !0.5! (YDT1)$) |- (YDT1);
        \draw[line] (XD0) -| ($(XD0) !0.5! (YDT1)$) |- (YDT1);
        \draw[line] (ZDT0) -| ($(ZDT0) !0.5! (YDT1)$) |- (YDT1);
        
        \draw[line] (YDT1)--(YDT2);
        \draw[line] (t1) -| ($(t1) !0.5! (YDT2)$) |- (YDT2);
        \draw[line] (DW1) -| ($(DW1) !0.5! (YDT2)$) |- (YDT2);
        \draw[line] (XD1) -| ($(XD1) !0.5! (YDT2)$) |- (YDT2);
        \draw[line] (ZDT1) -| ($(ZDT1) !0.5! (YDT2)$) |- (YDT2);
        
        \draw[line] (YDT2)--(YDTi);
        \draw[line] (t2) -| ($(t2) !0.5! (YDTi)$) |- (YDTi);
        \draw[line] (DW2) -| ($(DW2) !0.5! (YDTi)$) |- (YDTi);
        \draw[line] (XD2) -| ($(XD2) !0.5! (YDTi)$) |- (YDTi);
        \draw[line] (ZDT2) -| ($(ZDT2) !0.5! (YDTi)$) |- (YDTi);
        
        \draw[line] (YDTi)--(YDTNm1);
        \draw[line] (ti) -| ($(ti) !0.5! (YDTNm1)$) |- (YDTNm1);
        \draw[line] (DWi) -| ($(DWi) !0.5! (YDTNm1)$) |- (YDTNm1);
        \draw[line] (XDi) -| ($(XDi) !0.5! (YDTNm1)$) |- (YDTNm1);
        \draw[line] (ZDTi) -| ($(ZDTi) !0.5! (YDTNm1)$) |- (YDTNm1);
        
        \draw[line] (YDTNm1)--(YDTN);
        \draw[line] (tNm1) -| ($(tNm1) !0.5! (YDTN)$) |- (YDTN);
        \draw[line] (DWNm1) -| ($(DWNm1) !0.5! (YDTN)$) |- (YDTN);
        \draw[line] (XDNm1) -| ($(XDNm1) !0.5! (YDTN)$) |- (YDTN);
        \draw[line] (ZDTNm1) -| ($(ZDTNm1) !0.5! (YDTN)$) |- (YDTN);
        
        \draw[line] (YDTN)--(LDm);
        \draw[line] (YDN)--(LDm);
    \end{tikzpicture}  
    \caption{Architecture of the DBSDE scheme.}
    \label{fig1}
\end{figure}

We use the time discretization $\Delta$,  and for notational convenience we write $W_n = W_{t_n}$, $\Delta W_n = W_{n+1} - W_n$, $(X_n, Y_n, Z_n) = (X_{t_n}, Y_{t_n}, Z_{t_n})$,  $(X^{\Delta}_n, Y^{\Delta}_n, Z^{\Delta}_n)$ the approximations of $(X_n, Y_n, Z_n)$. By using the Euler-Maruyama method one obtains
\begin{align*}
    \begin{split}
         X^{\Delta}_{n+1} & = X^{\Delta}_n + a\left(t_n, X^{\Delta}_n\right) \Delta t + b\left(t_n, X^{\Delta}_n\right) \Delta W_n,\\
         Y^{\Delta}_{n+1} &= Y^{\Delta}_n - f\left(t_n, X^{\Delta}_n, Y^{\Delta}_n, Z^{\Delta}_n\right) \Delta t +  Z^{\Delta}_n \Delta W_n,
     \end{split}
\end{align*}
where $n = 0, \ldots, N-1.$ Since the Brownian motions are independent, $\Delta W_n \sim \mathcal{N}(\mathbf{0}_d,\,\Delta t \, \mathbf{I}_{d})$, with $\mathbf{0}_d \in \mathbb{R}^d$ a vector of zeros and $\mathbf{I}_{d} \in \mathbb{R}^{d \times d}$ the identity matrix. The discretized counterpart of~\eqref{eq4} is given as 
\begin{align}
    \begin{split}
        &\inf_{Y_0^{\Delta} \in \mathbb{L}^2_{\mathcal{F}_0}\left(\Omega;  \mathbb{R} \right),\, Z^{\Delta} \in \mathbb{H}^{\Delta, 2}\left(\{0, 1, \ldots, N-1\} \times \Omega; \mathbb{R}^{1 \times d} \right)} \mathbf{L}^{\Delta}\bigl(Y_0^{\Delta}, Z^{\Delta}\bigr),\\
         s.t. \quad & X^{\Delta}_{0} = x_0,\\
         & X^{\Delta}_{n+1} = X^{\Delta}_n + a\left(t_n, X^{\Delta}_n\right) \Delta t + b\left(t_n, X^{\Delta}_n\right) \Delta W_n,\\
         & Y^{\Delta}_{n+1} = Y^{\Delta}_n - f\left(t_n, X^{\Delta}_n, Y^{\Delta}_n, Z^{\Delta}_n\right) \Delta t +  Z^{\Delta}_n \Delta W_n,
    \end{split}
\label{eq5}
\end{align}
for $n = 0, 1, \ldots, N-1$, where 
$$\mathbf{L}^{\Delta}\bigl(Y_0^{\Delta}, Z^{\Delta}\bigr):= \mathbb{E}\bigl[|g(X^{\Delta}_N)-Y^{\Delta}_N|^2\bigr].$$ The authors in~\cite{weinan2017deep} considered DNNs, namely $\phi_0^y: \mathbb{R}^d \to \mathbb{R}$ and $\phi_n^z: \mathbb{R}^d \to \mathbb{R}^{1 \times d}$ for $n = 0, 1, \ldots, N-1$ to estimate the solution of~\eqref{eq5}. For a sample of size $m$ we have
\begin{align}
    \begin{split}
        &\min_{\theta \in \Theta} \mathbf{L}^{\Delta,m}\bigl(\phi_0^y(x_0; \theta_0^y), \phi^z(X^{\Delta}; \theta^z)\bigr),\\
        s.t. \quad & X^{\Delta}_{0} = x_0, \quad Y^{\Delta, \theta}_{0} = \phi_0^y(x_0; \theta_0^y),\\
         & X^{\Delta}_{n+1, j} = X^{\Delta}_{n, j}+ a\left(t_n, X^{\Delta}_{n, j}\right) \Delta t + b\left(t_n, X^{\Delta}_{n, j}\right) \sqrt{\Delta t} \mathcal{Z}_j,\\
         & Z^{\Delta, \theta}_{n, j} =  \phi_n^z(X^{\Delta}_{n, j}; \theta_n^z),\\
         & Y^{\Delta, \theta}_{n+1, j} = Y^{\Delta, \theta}_{n, j} - f\left(t_n, X^{\Delta}_{n, j}, Y^{\Delta, \theta}_{n, j}, Z^{\Delta, \theta}_{n, j}\right) \Delta t +  Z^{\Delta, \theta}_{n, j} \sqrt{\Delta t} \mathcal{Z}_j,
    \end{split}
\label{eq6}
\end{align}
for $n = 0, 1, \ldots, N-1$, $j=1, \ldots, m$, $\mathcal{Z}_j \sim \mathcal{N}(\mathbf{0}_d,\,\mathbf{I}_{d})$ and
$$\mathbf{L}^{\Delta,m}\bigl(\phi_0^y(x_0; \theta_0^y), \phi^z(X^{\Delta}; \theta^z)\bigr) := \frac{1}{m} \sum_{j=1}^m|g(X^{\Delta}_{N,j})-Y^{\Delta, \theta}_{N,j}|^2.$$
Note that $Y_0^{\Delta, \theta} =\theta_0^y \in \mathbb{R}$ and $Z_0^{\Delta, \theta}=\theta_0^z \in \mathbb{R}^{1 \times d}$ are considered as learnable parameters in~\cite{weinan2017deep}, which are initialized by sampling from uniform distributions. Furthermore, $N-1$ DNNs are employed to calculate $Z^{\Delta, \theta}_n$ for $n=1, 2, \ldots, N-1$. The architecture of the DBSDE scheme is displayed in Figure~\ref{fig1}. 

In the numerical section, we also consider the LaDBSDE scheme~\cite{kapllani2020deep} to demonstrate the applicability of our UQ model to other deep learning-based BSDE schemes. The LaDBSDE scheme addresses the issues encountered in the DBSDE scheme, such as convergence to an approximation far from the true solution or even divergence, especially for a complex solution structure and a long terminal time. In the LaDBSDE scheme, the BSDE problem is formulated as a global optimization problem with local loss functions at each time step. The process $Y$ is approximated using the same DNN $\phi^y$, while the process $Z$ is obtained through automatic differentiation due to~\eqref{eq3}. These approximations are performed by globally minimizing quadratic local loss functions defined at each time step, which always includes the terminal condition $Y_T$. The loss functions are obtained by iterating the Euler-Maruyama discretization of the integrals with the terminal condition $Y_T$. For further details, we refer to~\cite{kapllani2020deep}. In the following sections, we introduce all the sources of uncertainty in the DBSDE scheme and propose a UQ model to estimate the uncertainty. It is important to emphasize that our approach is applicable not only to the DBSDE scheme but also to the LaDBSDE scheme and potentially to other deep learning-based BSDE schemes as well.

\subsection{Sources of uncertainty in the DBSDE scheme}
\label{sec32}
For notational convenience, we consider
\begin{itemize}
    \item $\mathbf{Y}^{\star}:=\argmin_{\mathbf{Y} \in \mathbb{L}^2_{\mathcal{F}_0}\left(\Omega; \mathbb{R} \right) \times \mathbb{H}^2\left([0, T]  \times \Omega; \mathbb{R}^{1 \times d} \right)} \mathbf{L}(\mathbf{Y})$ the solution to~\eqref{eq4}, where $\mathbf{Y} = \left(Y_0, Z \right)$.
    \item $\mathbf{Y}^{\Delta,\star}:=\argmin_{\mathbf{Y}^{\Delta} \in \mathbb{L}^2_{\mathcal{F}_0}\left( \Omega; \mathbb{R} \right) \times \mathbb{H}^{\Delta, 2}\left( \{0, 1, \ldots, N-1\} \times \Omega; \mathbb{R}^{1 \times d} \right)} \mathbf{L}^{\Delta}(\mathbf{Y}^{\Delta})$ the solution to~\eqref{eq5}, where $\mathbf{Y}^{\Delta} = \left(Y^{\Delta}_0, Z^{\Delta} \right)$.
    \item $\theta^{\star}:=\argmin_{\theta \in \Theta} \mathbf{L}^{\Delta}(\mathbf{Y}^{\theta})$ the optimal parameters, where $\mathbf{Y}^{\theta} = \left(\phi_0^y(x_0; \theta_0^y), \phi^z(X^{\Delta}; \theta^z)\right)$. 
    \item $\theta^{m, \star}:=\argmin_{\theta \in \Theta} \mathbf{L}^{\Delta, m}(\mathbf{Y}^{\theta})$ the optimal parameters in~\eqref{eq6}.
    \item $\mathcal{A}: \mathbb{N} \to \Theta$ the optimization algorithm and $\hat{\theta}^m = \mathcal{A}(m)$ an estimate for the sample of size $m$.
    \item $\varphi^N: \mathbb{L}^2_{\mathcal{F}_0}\left( \Omega; \mathbb{R} \right) \times \mathbb{H}^{\Delta, 2}\left( \{0, \ldots, N-1\} \times \Omega; \mathbb{R}^{1 \times d} \right) \to  \mathbb{L}^2_{\mathcal{F}_0}\left( \Omega; \mathbb{R} \right) \times \mathbb{H}^2\left( [0, T] \times \Omega; \mathbb{R}^{1 \times d} \right)$ is an interpolation scheme. We choose $\varphi^N$ such that $\mathbf{L}^{\Delta} - \mathbf{L}\circ \varphi^N = 0$.
\end{itemize}
Then, constructing a telescopic sum we obtain the following error decomposition of the DBSDE scheme
\begin{equation*}
\begin{aligned}
    \mathbf{L}(\varphi^N(\mathbf{Y}^{\hat{\theta}^m})) &\leq \mathbf{L}(\mathbf{Y}^{\star})\\
      & + \mathbf{L}(\varphi^N(\mathbf{Y}^{\Delta,\star})) - \mathbf{L}(\mathbf{Y}^{\star})\\
      & + | \mathbf{L}^{\Delta}(\mathbf{Y}^{\Delta,\star}) - \mathbf{L}(\varphi^N(\mathbf{Y}^{\Delta,\star})) |\\
      & +  \mathbf{L}^{\Delta}(\mathbf{Y}^{\theta^{\star}})  - \mathbf{L}^{\Delta}(\mathbf{Y}^{\Delta,\star})\\
      & + \mathbf{L}^{\Delta}(\mathbf{Y}^{\theta^{m, \star}}) - \mathbf{L}^{\Delta}(\mathbf{Y}^{\theta^{\star}})\\
      & + | \mathbf{L}^{\Delta,m}(\mathbf{Y}^{\theta^{m, \star}}) - \mathbf{L}^{\Delta}(\mathbf{Y}^{\theta^{m, \star}}) |\\
      & +  \mathbf{L}^{\Delta,m}(\mathbf{Y}^{\hat{\theta}^m}) - \mathbf{L}^{\Delta,m}(\mathbf{Y}^{\theta^{m, \star}})\\
      & + | \mathbf{L}^{\Delta}(\mathbf{Y}^{\hat{\theta}^m}) - \mathbf{L}^{\Delta,m}(\mathbf{Y}^{\hat{\theta}^m}) |\\
      & + | \mathbf{L}(\varphi^N(\mathbf{Y}^{\hat{\theta}^m})) - \mathbf{L}^{\Delta}(\mathbf{Y}^{\hat{\theta}^m}) |.
\end{aligned}
\end{equation*}
Since $\mathbf{L}^{\Delta} - \mathbf{L}\circ \varphi^N = 0$, we have
\begin{align}
    \mathbf{L}(\varphi^N(\mathbf{Y}^{\hat{\theta}^m})) &\leq \mathbf{L}(\varphi^N(\mathbf{Y}^{\Delta,\star})) \label{eq7}\\
      & +  \mathbf{L}^{\Delta}(\mathbf{Y}^{\theta^{\star}})  - \mathbf{L}^{\Delta}(\mathbf{Y}^{\Delta,\star}) \label{eq8}\\
      & + \mathbf{L}^{\Delta}(\mathbf{Y}^{\theta^{m, \star}}) - \mathbf{L}^{\Delta}(\mathbf{Y}^{\theta^{\star}}) \label{eq9}\\
      & + | \mathbf{L}^{\Delta,m}(\mathbf{Y}^{\theta^{m, \star}}) - \mathbf{L}^{\Delta}(\mathbf{Y}^{\theta^{m, \star}}) | \label{eq10}\\
      & +  \mathbf{L}^{\Delta,m}(\mathbf{Y}^{\hat{\theta}^m}) - \mathbf{L}^{\Delta,m}(\mathbf{Y}^{\theta^{m, \star}}) \label{eq11}\\
      & + | \mathbf{L}^{\Delta}(\mathbf{Y}^{\hat{\theta}^m}) - \mathbf{L}^{\Delta,m}(\mathbf{Y}^{\hat{\theta}^m}) | \label{eq12}. 
\end{align}
In the decomposition above, each term corresponds to a specific source of uncertainty. The naming convention we adopt here is motivated by~\cite{hullermeier2021aleatoric}. The term~\eqref{eq7} captures the discretization error associated with the Euler-Maruyama method. Note that a DNN architecture has to be chosen when implementing the algorithm. Hence,~\eqref{eq8} represents the model or approximation error~\cite{beck2022full,jentzen2023strong}. Since the scheme optimizes the empirical loss,~\eqref{eq9} denotes the estimation error, as the empirical loss is only an estimate of the true loss. Selecting an SGD-type algorithm to optimize the DNN parameters introduces the optimization error, represented by~\eqref{eq11}. Finally,~\eqref{eq10} and \eqref{eq12} correspond to the sampling errors, which are insignificant compared to the other error sources. Identifying these errors is important to evaluate their practical impact on the uncertainty of the scheme. In the numerical section, we discuss such errors.

\subsection{The UQ model}
\label{sec33}
In practice, it is challenging to disentangle the sources of uncertainty in the DBSDE scheme, as we demonstrate in the numerical experiments. Consequently, quantifying the uncertainty of the DBSDE scheme becomes crucial for practical applications. In this section, we develop a UQ model to estimate the uncertainty. After applying the DBSDE scheme, we obtain the random variables $Y^{\Delta, \hat{\theta}^m}_0 \in \mathbb{R}$ and $Z^{\Delta, \hat{\theta}^m}_0 \in \mathbb{R}^{1 \times d}$ that approximate $Y_0$ and $Z_0$, respectively. To evaluate the quality of such approximations, one can use the expected squared error when the exact solution is known. This metric accounts for all the error sources in the DBSDE scheme and is calculated as 
\begin{equation*}
\epsilon^y:=\sqrt{\mathbb{E}\bigl[(Y_0^{\Delta, \hat{\theta}^m} - Y_0)^2\bigr]} \in \mathbb{R}^{+}, \quad \epsilon^{z_k}:=\sqrt{\mathbb{E}\bigl[(Z_0^{\Delta, \hat{\theta}^m, k} - Z_0^k)^2\bigr]}  \in \mathbb{R}^{+},
\end{equation*}
for $k = 1, \ldots, d$. However, $\left(Y_0, Z_0\right)$ is usually unknown, the STD of the approximate solutions
\begin{equation*}
   \sigma^y :=  \sqrt{\mathbb{E}\bigl[(Y_0^{\Delta, \hat{\theta}^m}-\mu^y)^2\bigr]}  \in \mathbb{R}^{+}, \quad \sigma^{z_k} :=  \sqrt{\mathbb{E}\bigl[(Z_0^{\Delta, \hat{\theta}^m, k}-\mu^{z_k})^2\bigr]}  \in \mathbb{R}^{+},
\end{equation*}
is often used, where $k = 1, \ldots, d,$ $\mu^y :=\mathbb{E}\bigl[Y_0^{\Delta, \hat{\theta}^m}\bigr]  \in \mathbb{R}$ and $\mu^{z_k} := \mathbb{E}\bigl[Z_0^{\Delta, \hat{\theta}^m, k}\bigr]  \in \mathbb{R}$. To compute the STD (and the expected squared error when the exact solution $\left(Y_0, Z_0\right)$ is available), $Q$ runs of the DBSDE algorithm must be done. The STD and the expected squared error are used as the benchmark in our experiments. To this end, we have the root mean squared error (RMSE) and the ensemble (biased sample) STD as
\begin{equation*}
    \tilde{\epsilon}^y := \sqrt{\frac{1}{Q}\sum_{q=1}^{Q} \left( Y_{0,q}^{\Delta, \hat{\theta}^m} - Y_0\right)^2 }, \quad \tilde{\sigma}^y := \sqrt{\frac{1}{Q}\sum_{q=1}^{Q} \left( Y_{0,q}^{\Delta, \hat{\theta}^m} - \tilde{\mu}^y\right)^2 },
\end{equation*}
for $Y_0$ and 
\begin{equation*}
    \tilde{\epsilon}^{z_k} := \sqrt{\frac{1}{Q}\sum_{q=1}^{Q} \left( Z_{0,q}^{\Delta, \hat{\theta}^m, k} - Z_0^k\right)^2 }, \quad \tilde{\sigma}^{z_k} := \sqrt{\frac{1}{Q}\sum_{q=1}^{Q} \left( Z_{0,q}^{\Delta, \hat{\theta}^m, k} - \tilde{\mu}^{z_k}\right)^2 },
\end{equation*}
for $Z_0$, $k = 1, \ldots, d$, where $\left(Y_{0,q}^{\Delta, \hat{\theta}^m}, Z_{0,q}^{\Delta, \hat{\theta}^m}\right)$ represents the approximated solutions from the $q$-th run of the algorithm, $\tilde{\mu}^y:=\frac{1}{Q}\sum_{q=1}^{Q}Y_{0,q}^{\Delta, \hat{\theta}^m}$ and $\tilde{\mu}^{z_k}:=\frac{1}{Q}\sum_{q=1}^{Q}Z_{0,q}^{\Delta, \hat{\theta}^m,k}$ are the ensemble (sample) means of the approximate solution. Note that one can use the ensemble unbiased STD. However, the estimate of the STD from the UQ model is biased (as a maximum likelihood estimate). Therefore, for the purpose of comparison in numerical experiments, it is more convenient to use the ensemble biased STD. Usually, $Q=10$ is used, which is computationally expensive in high dimensions. Therefore, we propose a UQ model to estimate the STDs for $\left(Y_0^{\Delta, \hat{\theta}^m}, Z_0^{\Delta, \hat{\theta}^m}\right)$ by using only $Q=1$ run of the algorithm. The model is based on an approach commonly used to quantify uncertainty in heteroscedastic nonlinear regression. We make the assumption that the errors of the DBSDE scheme follow a normal distribution with zero mean and the STD depending on the parameter set of the discretized BSDE (such as $T$, $x_0$, $\Delta t$, etc.). This assumption aligns with the standard practice in heteroscedastic regression. To train the UQ model, we construct a dataset of length $M$ with i.i.d samples $\mathcal{D} = \{ \mathbf{x}_i, \mathbf{y}_i, \mathbf{z}_i \}_{i=1}^{M}$. Here, $\mathbf{x}_i \in \mathbb{R}^{\nf}$ represents $\nf$ parameters of the discretized BSDE, for example, $\mathbf{x}_i :=(x_{0,i}, T_i, \Delta t_i)$, and $\left(\mathbf{y}_i, \mathbf{z}_i\right) := \left(Y_0^{\Delta, \hat{\theta}^m}(\mathbf{x}_i), Z_0^{\Delta, \hat{\theta}^m}(\mathbf{x}_i)\right) \in \mathbb{R} \times \mathbb{R}^{1 \times d}$ are the approximations of $(Y_0(\mathbf{x}_i), Z_0(\mathbf{x}_i))$ obtained from the DBSDE scheme using parameter set $\mathbf{x}_i$. We use uniform distributions to select the parameter set $\mathbf{x}_i$
$$
x_{0, i} \sim \mathcal{U}[x_0^{min}, \, x_0^{max}], \quad T_{i} \sim \mathcal{U}[T^{min}, \, T^{max}],
$$
where $(x_0^{min}, x_0^{max})$ and $(T^{min}, T^{max})$ are the boundaries of the corresponding uniform distributions for $x_0$ and $T$. The dataset $\mathcal{D}$ is generated using Algorithm~\ref{alg1}.
\begin{algorithm}
\caption{Algorithm generating dataset $\mathcal{D}$ for UQ model}\label{alg1}
\KwIn{$\left(N, M, d, x_0^{min}, x_0^{max}, T^{min}, T^{max}, Y_0^{min}, Y_0^{max}\right)$ - problem related parameters} 
\KwIn{$\left(a, b, f, g\right)$ - functions of BSDE system}
\KwIn{$\left(\alpha, \Kf, L, \eta, \varrho, m\right)$ - DNN hyperparameters in DBSDE scheme}
\KwOut{$\mathcal{D} = \{ \mathbf{x}_i, \mathbf{y}_i, \mathbf{z}_i \}_{i=1}^{M}$ - Dataset for UQ model}
\For{$i= 1: M$}{
    $x_{0, i} \sim \mathcal{U}[x_0^{min}, \, x_0^{max}]$\\
    $T_{i} \sim \mathcal{U}[T^{min}, \, T^{max}]$\\
    $\Delta t_i = \frac{T_i}{N}$\\
    $\mathbf{x}_i =(x_{0,i}, T_i, \Delta t_i)$\\
    \For{$n= 0: N$}{
        $t_{n} = n \,\Delta t_i$\\
    }
    \textbf{\underline{Initialize parameter set $\theta$}}\\
    $\hat{\theta}^{y, m, 0}_0 = Y_0^{\Delta, \hat{\theta}^{m, 0}}(\mathbf{x}_i) \sim \mathcal{U}[Y_0^{min}, \, Y_0^{max}]$\\
    $\hat{\theta}^{z, m, 0}_0 =Z_0^{\Delta, \hat{\theta}^{m, 0}}(\mathbf{x}_i) \sim \mathcal{U}[-\mathbf{1}_d,\, \mathbf{1}_d]$ - $\mathbf{1}_d$ vector of all ones\\
    $\left(\hat{\theta}^{z, m, 0}_1, \ldots,  \hat{\theta}^{z, m, 0}_{N-1}\right)$- Xavier normal initializer~\cite{glorot2010understanding}\\    
    $\hat{\theta}^{m, 0} = \left( \hat{\theta}^{y, m, 0}_0, \hat{\theta}^{z, m, 0}_0, \hat{\theta}^{z, m, 0}_1, \ldots, \hat{\theta}^{z, m, 0}_{N-1}\right)$\\
    \textbf{\underline{Optimization or training part}}\\
    \For{$\kappa= 1: \Kf$}{
        \For{$j = 1: m$}{
            $X^{\Delta}_{0,j} = x_{0,i}$\\
            \For{$n = 0: N-1$}{
                \textbf{\underline{Euler-Maruyama for the forward SDE}}\\
                $\Delta W_{n,j} \sim \mathcal{N}(\mathbf{0}_d, \Delta t_i \,\mathbf{I}_d)$\\
                $X^{\Delta}_{n+1,j} = X^{\Delta}_{n,j} + a\left(t_n, X^{\Delta}_{n, j}\right) \Delta t_i + b\left(t_n, X^{\Delta}_{n,j}\right) \Delta W_{n,j}$\\
                \textbf{\underline{Use DNN with $\left(L, \eta, \varrho\right)$ for $Z$ and Euler-Maruyama for $Y$}}\\
                $d_0 = d_1 = d$\\
                \eIf{$n < N-1$}{
                    $Z^{\Delta, \hat{\theta}^{m, \kappa-1}}_{n+1, j} =  \phi^z_{n+1}(X^{\Delta}_{n+1,j};\hat{\theta}^{z, m, \kappa-1}_{n+1})$\\
                     $Y^{\Delta, \hat{\theta}^{m, \kappa-1}}_{n+1, j} = Y^{\Delta, \hat{\theta}^{m, \kappa-1}}_{n, j} - f(t_n, X^{\Delta}_{n,j}, Y^{\Delta, \hat{\theta}^{m, \kappa-1}}_{n, j}, Z^{\Delta, \hat{\theta}^{m, \kappa-1}}_{n, j})\Delta t_i$\\
                    $\quad \quad \quad \quad \quad \,\,+ Z^{\Delta, \hat{\theta}^{m, \kappa-1}}_{n, j} \Delta W_{n, j}$\\                
                }{
                  $Y^{\Delta, \hat{\theta}^{m, \kappa-1}}_{n+1, j} = Y^{\Delta, \hat{\theta}^{m, \kappa-1}}_{n, j} - f(t_n, X^{\Delta}_{n,j}, Y^{\Delta, \hat{\theta}^{m, \kappa-1}}_{n, j}, Z^{\Delta, \hat{\theta}^{m, \kappa-1}}_{n, j})\Delta t_i$\\
                    $\quad \quad \quad \quad \quad \,\,+ Z^{\Delta, \hat{\theta}^{m, \kappa-1}}_{n, j} \Delta W_{n, j}$\\ 
                }
            }
        }
        $\mathbf{L}^{\Delta, m}(\hat{\theta}^{m, \kappa-1}) = \frac{1}{m}\sum_{j=1}^m | g(X^{\Delta}_{N,j}) -  Y^{\Delta, \hat{\theta}^{m, \kappa-1}}_{N, j} |^2$\\      
         \textbf{\underline{Adam optimization step}}\\
        $\hat{\theta}^{m,\kappa}$ - trained parameters with Adam optimizer~\cite{kingma2014adam}, learning rate $\alpha$\\        
    }
  $\hat{\theta}^{m} = \hat{\theta}^{m,\Kf}$ - final estimated parameters after $\Kf$ optimization steps\\
  $\left(\mathbf{y}_i, \mathbf{z}_i\right) = \left(Y_0^{\Delta, \hat{\theta}^{m}}(\mathbf{x}_i), Z_0^{\Delta, \hat{\theta}^{m}}(\mathbf{x}_i)\right)$\\
}
\end{algorithm}
Note that one can use the entire dataset to build a learning algorithm that considers $Y_0$ and $Z_0$ as pairs (the BSDE solution at $t_0$ is the pair $\left(Y_0, Z_0\right)$). However, this approach may introduce increased complexity for the learning algorithm, mainly because the magnitudes of the solutions for $\mathbf{y}$ and $\mathbf{z}$ over different parameter sets $\mathbf{x}$ can differ significantly. Additionally, assumptions regarding their correlation might be necessary. Hence, we divide the dataset $\mathcal{D}$ into two datasets, $\mathcal{D}^y = \{ \mathbf{x}_i, \mathbf{y}_i \}_{i=1}^{M}$ and $\mathcal{D}^z = \{ \mathbf{x}_i, \mathbf{z}_i \}_{i=1}^{M}$, and develop two different learning algorithms. Given the input feature $\mathbf{x}$, we use one DNN to model the probabilistic predictive distribution $p^{\omega}(\mathbf{y}|\mathbf{x})$ and another for $p^{\psi}(\mathbf{z}|\mathbf{x})$, where $\omega$ and $\psi$ are the parameters of the corresponding DNNs. More precisely, we treat each observed value as a sample from a Gaussian distribution (multivariate Gaussian for $Z_0$), with the mean and STD as functions of the parameter set, namely
\begin{equation}
    \mathbf{y}_i \sim \mathcal{N}(\mu^y\left(\mathbf{x}_i\right),\, (\sigma^y)^2\left(\mathbf{x}_i\right)), \quad \mathbf{z}_i \sim \mathcal{N}(\mu^z\left(\mathbf{x}_i\right),\, \mathbf{I}_d\,(\sigma^z)^2\left(\mathbf{x}_i\right)),
    \label{eq13}
\end{equation}
and allow the networks to calculate their estimates as $\left(\hat{\mu}^{y,\omega}(\mathbf{x}_i), \hat{\sigma}^{y,\omega}(\mathbf{x}_i)\right) \in \mathbb{R} \times \mathbb{R}^{+}$ and $\left(\hat{\mu}^{z,\psi}(\mathbf{x}_i), \hat{\sigma}^{z,\psi}(\mathbf{x}_i)\right) \in \mathbb{R}^d \times \mathbb{R}^{d,+}$. These calculations are performed by minimizing the negative log-likelihood 
\begin{equation*}
    \mathbf{L}^{\omega, M}(\mathcal{D}^y) := -\log \left( p^{\omega}(\mathbf{y}|\mathbf{x})\right) = \frac{1}{M}\sum_{i=1}^M \Biggl(  \log \left( \hat{\sigma}^{y,\omega}(\mathbf{x}_i) \right) +  \frac{1}{2}\frac{\left(\mathbf{y}_i - \hat{\mu}^{y,\omega}(\mathbf{x}_i)\right)^2}{(\hat{\sigma}^{y,\omega})^2(\mathbf{x}_i)} \Biggr) + c^y,
\end{equation*}
for $Y_0$ and assuming that the covariance matrix of $Z^{\Delta, \hat{\theta}^m}_0$ is diagonal, then
\begin{equation*}
    \begin{split}
        \mathbf{L}^{\psi,M}(\mathcal{D}^z) := -\log \left( p^{\psi}(\mathbf{z}|\mathbf{x})\right) & = \frac{1}{M}\sum_{i=1}^M \Biggl( \sum_{k=1}^{d}\log \left( \hat{\sigma}^{z_k,\psi}(\mathbf{x}_i) \right) \Biggr. \\
        & \quad + \Biggr. \frac{1}{2} \left(\mathbf{z}_i - \hat{\mu}^{z,\psi}(\mathbf{x}_i)\right)^{\top} \left((\hat{\sigma}^{z,\psi})^2(\mathbf{x}_i)\right)^{-1} \left(\mathbf{z}_i - \hat{\mu}^{z,\psi}(\mathbf{x}_i)\right) \Biggr) + c^z,
    \end{split}
\end{equation*}
for $Z_0$, where $c^y>0, c^z>0$ are constants. We use Algorithm~\ref{alg2} and~\ref{alg3} to estimate the parameters in~\eqref{eq13} for $Y_0$ and $Z_0$, respectively.
\begin{algorithm}
\caption{Algorithm estimating the parameters of~\eqref{eq13} for $Y_0$}\label{alg2}
\KwIn{$\left(M, \nf, M^{valid}, M^{test}, \mathcal{D}^y\right)$ - parameters and dataset for UQ model}
\KwIn{$\left(L^y, \eta^y, \varrho^y, \alpha^y, B^y, ep^y, \lambda^y\right)$ - DNN hyperparameters}
\KwOut{$\left( \hat{\mu}^{y,\hat{\omega}^M}(\mathbf{x}), \hat{\sigma}^{y,\hat{\omega}^M}(\mathbf{x})\right)$ - estimates of UQ model for $Y_0$}
\textbf{\underline{Split $\mathcal{D}^y$ into training, validation and testing samples}}\\   
$M^{train} = M - M^{valid} - M^{test}$\\
$(\mathbf{x}^{train}, \mathbf{y}^{train})$\\
$(\mathbf{x}^{valid}, \mathbf{y}^{valid})$\\
$(\mathbf{x}^{test}, \mathbf{y}^{test})$\\
\textbf{\underline{Normalize input data $\mathbf{x}$ based on training statistics}}\\    
$(\mathbf{x}^{train, nr}, \mathbf{x}^{valid, nr}, \mathbf{x}^{test, nr})$\\
$d_0 = \nf$\\    
$d_1^y = 2$\\     
\textbf{\underline{Initialize parameters $\omega$}}\\
$\hat{\omega}^{M, 0}$- Xavier normal initializer~\cite{glorot2010understanding}\\ 
$\kappa = 0$\\
\For{$e=1:ep^y$}{
     $\kappa = \kappa+1$\\
    \For{$I=1:\frac{M^{train}}{B^y}$}{
        \textbf{\underline{Batch data}}\\
        $\mathcal{D}^{y,train,nr} = \{\mathbf{x}_i^{train,nr}, \mathbf{y}_i^{train}\}_{i=(I-1)\,B^y+1}^{I\,B^y}$\\        \textbf{\underline{Use DNN with $\left(L^y, \eta^y, \varrho^y\right)$ to estimate parameters in~\eqref{eq13}} for $Y_0$}\\
        $\left( \hat{\mu}^{y,\hat{\omega}^{M, \kappa-1}}(\mathbf{x}), \hat{\sigma}^{y,\hat{\omega}^{M, \kappa-1}}(\mathbf{x})\right) = \phi^y\left( \{\mathbf{x}_i^{train,nr}\}_{i=(I-1)\,B^y+1}^{I\,B^y}; \hat{\omega}^{M, \kappa-1}\right)$\\
        \textbf{\underline{Calculate loss including L2 regularization}}\\
         $\mathbf{L}^{\hat{\omega}^{M, \kappa-1}, M}(\mathcal{D}^{y,train,nr}) = \frac{1}{B^y}\sum_{i=(I-1)\,B^y+1}^{I\,B^y} \Biggl( 
         \log \left( \hat{\sigma}^{y,\hat{\omega}^{M, \kappa-1}}(\mathbf{x}_i^{train,nr}) \right) \Biggr.$ \\
         $\qquad \qquad \qquad \qquad \qquad \,\, + \Biggr.  \frac{1}{2}\frac{\left(\mathbf{y}_i^{train,nr} - \hat{\mu}^{y,\hat{\psi}^{M, \kappa-1}}(\mathbf{x}_i^{train,nr})\right)^2}{(\hat{\sigma}^{y,\hat{\omega}^{M, \kappa-1}})^2(\mathbf{x}_i^{train,nr})}\Biggr)$\\
         $\qquad \qquad \qquad \qquad \qquad \,\, + \lambda^y \sum_{l=1}^{L^y+1} \left(\hat{\omega}^{M, \kappa-1}(l)\right)^2$\\
         \textbf{\underline{Adam optimization step}}\\
         $\hat{\omega}^{M, \kappa}$ - trained parameters with Adam optimizer~\cite{kingma2014adam}, learning rate $\alpha^y$\\
    }
}
$\hat{\omega}^{M} = \hat{\omega}^{M, \kappa}$ - final estimated parameters of DNN after $ep^y$ epochs, each with $\frac{M^{train}}{B^y}$ number of batches\\
$\left( \hat{\mu}^{y,\hat{\omega}^M}(\mathbf{x}), \hat{\sigma}^{y,\hat{\omega}^M}(\mathbf{x})\right)$ - estimated parameters of~\eqref{eq13} for $Y_0$, $\mathbf{x}$ training, validation or testing sample
\end{algorithm}
\begin{algorithm}
\caption{Algorithm estimating the parameters of~\eqref{eq13} for $Z_0$}\label{alg3}
\KwIn{$\left(M, \nf, M^{valid}, M^{test}, \mathcal{D}^z\right)$ - parameters and dataset for UQ model}
\KwIn{$\left(L^z, \eta^z, \varrho^z, \alpha^z, B^z, ep^z, \lambda^z\right)$ - DNN hyperparameters}
\KwOut{$\left( \hat{\mu}^{z,\hat{\psi}^M}(\mathbf{x}), \hat{\sigma}^{z,\hat{\psi}^M}(\mathbf{x})\right)$ - estimates of UQ model for $Z_0$}
\textbf{\underline{Split $\mathcal{D}^z$ into training, validation and testing samples}}\\   
$M^{train} = M - M^{valid} - M^{test}$\\
$(\mathbf{x}^{train}, \mathbf{z}^{train})$\\
$(\mathbf{x}^{valid}, \mathbf{z}^{valid})$\\
$(\mathbf{x}^{test}, \mathbf{z}^{test})$\\
\textbf{\underline{Normalize input data $\mathbf{x}$ based on training statistics}}\\    
$(\mathbf{x}^{train, nr}, \mathbf{x}^{valid, nr}, \mathbf{x}^{test, nr})$\\
$d_0 = \nf$\\    
$d_1^z = 2d$\\     
\textbf{\underline{Initialize parameters $\psi$}}\\
$\hat{\psi}^{M, 0}$- Xavier normal initializer~\cite{glorot2010understanding}\\ 
$\kappa=0$\\
\For{$e=1:ep^z$}{
    $\kappa = \kappa+1$\\
    \For{$I=1:\frac{M^{train}}{B^z}$}{
        \textbf{\underline{Batch data}}\\
        $\mathcal{D}^{z,train,nr} = \{\mathbf{x}_i^{train,nr}, \mathbf{z}_i^{train}\}_{i=(I-1)\,B^z+1}^{I\,B^z}$\\        \textbf{\underline{Use DNN with $\left(L^z, \eta^z, \varrho^z\right)$ to estimate parameters in~\eqref{eq13}} for $Z_0$}\\
        $\left( \hat{\mu}^{z,\hat{\psi}^{M, \kappa-1}}(\mathbf{x}), \hat{\sigma}^{z,\hat{\psi}^{M, \kappa-1}}(\mathbf{x})\right) = \phi^z\left( \{\mathbf{x}_i^{train,nr}\}_{i=(I-1)\,B^z+1}^{I\,B^z}; \hat{\psi}^{M, \kappa-1}\right)$\\
        \textbf{\underline{Calculate loss including L2 regularization}}\\
         $\mathbf{L}^{\hat{\psi}^{M, \kappa-1}, M}(\mathcal{D}^{z,train,nr}) = \frac{1}{B^z}\sum_{i=(I-1)\,B^z+1}^{I\,B^z} \Biggl( \sum_{k=1}^{d}\log \left( \hat{\sigma}^{z_k,\hat{\psi}^{M, \kappa-1}}(\mathbf{x}_i^{train,nr}) \right) \Biggr.$ \\
         $\qquad \qquad \qquad \qquad \qquad \,\, + \Biggr. \frac{1}{2} \left(\mathbf{z}_i^{train,nr} - \hat{\mu}^{z,\hat{\psi}^{M, \kappa-1}}(\mathbf{x}_i^{train,nr})\right)^{\top}\Biggr.$ \\
         $\qquad \qquad \qquad \qquad \qquad \,\, \Biggr. \left((\hat{\sigma}^{z,\hat{\psi}^{M, \kappa-1}})^2(\mathbf{x}_i^{train,nr})\right)^{-1} \left(\mathbf{z}_i^{train,nr} - \hat{\mu}^{z,\hat{\psi}^{M, \kappa-1}}(\mathbf{x}_i^{train,nr})\right) \Biggr)$\\
         $\qquad \qquad \qquad \qquad \qquad \,\, + \lambda^z \sum_{l=1}^{L^z+1} \left(\hat{\psi}^{M, \kappa-1}(l)\right)^2$\\
         \textbf{\underline{Adam optimization step}}\\
         $\hat{\psi}^{M, \kappa}$ - trained parameters with Adam optimizer~\cite{kingma2014adam}, learning rate $\alpha^z$\\
    }
}
$\hat{\psi}^{M} = \hat{\psi}^{M, \kappa}$ - final estimated parameters of DNN after $ep^z$ epochs, each with $\frac{M^{train}}{B^z}$ number of batches\\
$\left( \hat{\mu}^{z,\hat{\psi}^M}(\mathbf{x}), \hat{\sigma}^{z,\hat{\psi}^M}(\mathbf{x})\right)$ - estimated parameters of~\eqref{eq13} for $Z_0$, $\mathbf{x}$ training, validation or testing sample
\end{algorithm}

\section{Numerical results}
\label{sec4}
In this section, we take the DBSDE scheme as an example to illustrate the impact of different sources of uncertainty in the scheme and apply our UQ model to both the DBSDE and LaDBSDE schemes. All the experiments were run in PYTHON using TensorFlow on the PLEIADES cluster, which consists of 268 workernodes. Each workernode has 2 sockets with an AMD EPYC 7452 32-Core processor (256GB of memory). For more information, see~PLEIADES documentation\footnote{\href{https://pleiadesbuw.github.io/PleiadesUserDocumentation/}{https://pleiadesbuw.github.io/PleiadesUserDocumentation/}}. 

\subsection{Experimental setup}
\label{subsec41}
To efficiently generate the dataset for the UQ model, we consider a straightforward parallelization of Algorithm~\ref{alg1}, namely one core is used to simulate the DBSDE scheme for one parameter
set. Hence, multiple cores provide a parallel generation of the dataset $\mathcal{D}$. The implementation of Algorithm~\ref{alg1} follows the hyperparameters considered in~\cite{weinan2017deep} for the DBSDE scheme. Each of the DNNs consists of $L=2$ hidden layers with $\eta=10+d$ neurons per hidden layer, and the rectifier function (ReLU) $\varrho(x)=\max(0, x) \in [0, \infty)$ is used as the activation function. Batch normalization is applied right after each matrix multiplication and before activation. Furthermore, all network parameters are initialized using a normal distribution without any pre-training. The Adam optimizer with learning rate $\alpha$ and $\Kf$ optimization steps is used as an SGD-type algorithm. For the implementation of Algorithm~\ref{alg2}, we first split the dataset $\mathcal{D}^y$ into training, validation, and testing samples, whose sizes are denoted by $M^{train}$, $M^{valid}$ and $M^{test}$, respectively. Note that the validation set is used to validate the performance of our UQ model when tuning its hyperparameters. The input layer of the DNN has $\nf$ neurons, and the output layer has $2$ neurons. The first neuron in the output layer estimates the mean of the approximate solution $\mu^y(\mathbf{x})$, and the second neuron in the output layer estimates the STD of the approximate solution $\sigma^y(\mathbf{x})$, where the softplus activation function $\varrho(x)=\ln(1 + e^{x}) \in (0, \infty)$ is applied to obtain positive estimates. The ReLU activation function is used for the $L^y$ hidden layers. Note that it is appropriate to choose $\eta^y > d_1^y$. The input data $\mathbf{x}$ is normalized based on the training data. We use the Adam optimizer with a batch size $B^y$, L2 regularization with parameter $\lambda^y$, and a specified number of epochs $ep^y$. The hyperparameters are set in a similar fashion for the implementation of Algorithm~\ref{alg3}.

To visually and quantitatively compare our estimates of the mean and STD of the approximate solution to benchmark values, such as the exact solution, RMSE, ensemble mean, and ensemble STD, we consider both linear and nonlinear BSDEs with available analytical solutions. We take the Black-Scholes BSDE as a linear $1$-dimensional example, which is used for pricing the European options.
\begin{exmp}
The Black-Scholes BSDE reads~\cite{zhao2010stable}
\begin{equation*}
    \begin{split}
    \left\{
        \begin{array}{rcl}
             dS_t &=&  a S_t\,dt + b S_t \,dW_t, \quad S_0 = s_0,\\  
   		   	   -dY_t &=& - \left( R Y_t + \left(  a - R + \delta \right)\frac{Z_t}{b}\right)\,dt- Z_t \,dW_t,\\ 
   		   	   Y_T &=& \left(S_T-K\right)^+.
        \end{array}
    \right. \\ 
    \end{split}
\end{equation*}\\[-7.5ex]
\label{ex1}
\end{exmp}
Note that $a$ represents the expected return of the stock $S_t$, $b$ denotes the volatility of the stock returns, $\delta$ is the dividend rate the stock pays, and $S_0$ the price of the stock at $t =0$. Moreover, $T$ denotes the maturity of the option contract, while $K$ represents the contract's strike price. Finally, $R$ corresponds to the risk-free interest rate. The analytic solution (the option price $Y_t$ and its delta hedging strategy $Z_t$) is given by
\begin{equation*}
 \begin{split}
    \left\{
        \begin{array}{rcl}
            	Y_t &=& S_t \exp\left(-\delta \left(T-t\right)\right) \Phi \left(d_1\right)-K \exp\left(-R\left(T-t\right)\right) \Phi \left(d_2\right),\\
   		   	   Z_t &=& S_t \exp\left(-\delta \left(T-t\right)\right) \Phi \left(d_1\right)b, \\
   		   	   d_{1/2} &=& \frac{\ln\left(\frac{S_t}{K}\right) + \left( R-\delta \pm \frac{b^2}{2} \right) \left(T-t\right)}{b \sqrt{T-t}},
        \end{array}
    \right. \\
    \end{split}
\end{equation*}
where $\Phi \left(\cdot\right)$ is the standard normal cumulative distribution function. For the nonlinear case, we consider the nonlinear high-dimensional Burgers type BSDE.
\begin{exmp}
The Burgers type BSDE reads~\cite{weinan2017deep}
\begin{equation*}
   		  \begin{split}
        \left\{
            \begin{array}{rcl}
                dX_t & = & b \, dW_t, \quad X_0 = 0,\\ 
                -dY_t & = & \left( \frac{b}{d} Y_t - \frac{2d+b^2}{2bd} \right) \left( \sum_{k=1}^{d} Z_t^k\right)\,dt -Z_t \,dW_t,\\  
   	   	    	Y_T & = & \frac{\exp\left( T +\frac{1}{d} \sum_{k=1}^{d} X_T^k\right)}{1+\exp\left( T +\frac{1}{d} \sum_{k=1}^{d} X_T^k\right)},
            \end{array}
        \right. 
    \end{split}
\end{equation*}\\[-7.5ex]
\label{ex2}
\end{exmp}
where, $W_t = (W_t^1, W_t^2, \ldots, W_t^d)^{\top}$, $X_t = (X_t^1, X_t^2, \ldots, X_t^d)^{\top}$ and $Z_t = (Z_t^1, Z_t^2, \ldots, Z_t^d)$.
The analytic solution is given by
\begin{equation*}
    \begin{split}
        \left\{
            \begin{array}{rcl}
                Y_t & = & \frac{\exp\left( t +\frac{1}{d} \sum_{k=1}^{d} X_t^k\right)}{1+\exp\left( t +\frac{1}{d} \sum_{k=1}^{d} X_t^k\right)},\\
   	   	    	Z_t & = & \frac{b}{d} \frac{\exp\left( t +\frac{1}{d} \sum_{k=1}^{d} X_t^k\right)}{\left(1+\exp\left( t +\frac{1}{d} \sum_{k=1}^{d} X_t^k\right)\right)^2}\mathbf{1}_{d}.
            \end{array}
        \right. 
    \end{split}
\end{equation*}
At $t_0$ we have that $\left(Y_0, Z_0 \right) = \left(0.5, \frac{b}{4d}\mathbf{1}_{d}\right)$.

The following experiments are organized as follows. Firstly, we illustrate the impact of the sources of uncertainty in the DBSDE scheme for both examples by visualizing the effect of different errors on the uncertainty. Secondly, we assess the performance of our UQ model for the mean and STD of the approximate solution by comparing them with the benchmark values. Additionally, we determine the number of runs of the DBSDE algorithm for which the ensemble STD is comparable to the estimated STD. Furthermore, the computational cost of generating the training data for the UQ model is evaluated. To demonstrate the applicability of the UQ model to other deep learning-based BSDE schemes, we apply it to the LaDBSDE scheme. Finally, we show the practical implications of our UQ model.

\subsection{The impact of the sources of uncertainty in the DBSDE scheme}
\label{subsec42}
To demonstrate the impact of the sources of uncertainty in the DBSDE scheme, we fix the parameter set of the BSDE and vary the hyperparameters of the DBSDE scheme that affect the corresponding source of uncertainty. Initially, we focus on the estimation and optimization errors. By using a high number of optimization steps $\Kf$ and different learning rate approaches, the effect of these errors on the uncertainty of the DBSDE scheme is analyzed. Afterward, we investigate the impact of the model error by increasing the number of hidden neurons $\eta$. Lastly, the number of discretization points $N$ is varied in order to study how the discretization error contributes to the uncertainty of the DBSDE scheme.

For the parameter values $T = 1, K = 100, S_0 = 100, a = 0.05, b = 0.2, R = 0.03$ and $\delta = 0$ in Example~\ref{ex1}, the exact solution is $(Y_0, Z_0) = (9.4134, 11.9741)$. In Example~\ref{ex2}, we chose $d=50$ and fix $b=25, T = 0.25$. The exact solution is  $\left( Y_0, Z_0\right) = \left(0.5, 0.125 \mathbf{1}_{50}\right)$. For the DBSDE algorithm, we start with the following hyperparameters: a constant learning rate approach (C-LR) with $\Kf = 60000$ optimization steps of the Adam algorithm, a learning rate $\alpha = \num{1e-2}$, and a batch size of $m=128$. To analyze the effect of estimation and optimization errors on the RMSE, we plot the RMSE values in Figures~\ref{fig2} and~\ref{fig3} for Examples~\ref{ex1} and~\ref{ex2}, respectively, for increasing values of $\Kf$. The RMSE values of $Z_0$ are plotted only for the first component in Example~\ref{ex2} as it is similar for the other components in our experiments. We use $N = 32$ discretization points and $Q = 10$ runs of the DBSDE algorithm. Note that each run of the DBSDE algorithm involves a different seed for generating the dataset and different initialization values of DNN parameters.
\begin{figure}[htb!]
        \centering
	\begin{subfigure}[h!]{0.48\linewidth}
            \begin{tikzpicture} 
            \begin{axis}[
                xmin = 0, xmax = 60000,
                ymin = 1e-2, ymax = 10,
                xtick distance = 10000,
                ymode=log,
                grid = both,
                width = \textwidth,
                height = 0.8\textwidth,
                major grid style = {lightgray},
                minor grid style = {lightgray!25},
                xlabel = {$\Kf$},
                ylabel = {$\tilde{\epsilon}^y$},]
            \addplot[
                smooth,
                ultra thick,
                blue,
                solid
            ] file {"Figures/Example1/Base/RMSE_Y_0.dat"};
            \end{axis}
            \end{tikzpicture}
		\caption{RMSE values for $Y_0$.}
		\label{fig2a}
	\end{subfigure}
	\begin{subfigure}[h!]{0.48\linewidth}
            \begin{tikzpicture} 
            \begin{axis}[
                xmin = 0, xmax = 60000,
                ymin = 1e-2, ymax = 10,
                xtick distance = 10000,
                ymode=log,
                grid = both,
                width = \textwidth,
                height = 0.8\textwidth,
                major grid style = {lightgray},
                minor grid style = {lightgray!25},
                xlabel = {$\Kf$},
                ylabel = {$\tilde{\epsilon}^z$},]
            \addplot[
                smooth,
                ultra thick,
                blue,
                solid
            ] file {"Figures/Example1/Base/RMSE_Z_0.dat"};
            \end{axis}
            \end{tikzpicture}
		\caption{RMSE values for $Z_0$.}
		\label{fig2b}
	\end{subfigure}
	\caption{RMSE values are plotted for Example~\ref{ex1} for increasing $\Kf,$ where $T = 1, K = 100, S_0 = 100, a = 0.05, b = 0.2, R = 0.03$ and $\delta = 0$.}
	\label{fig2}
\end{figure}
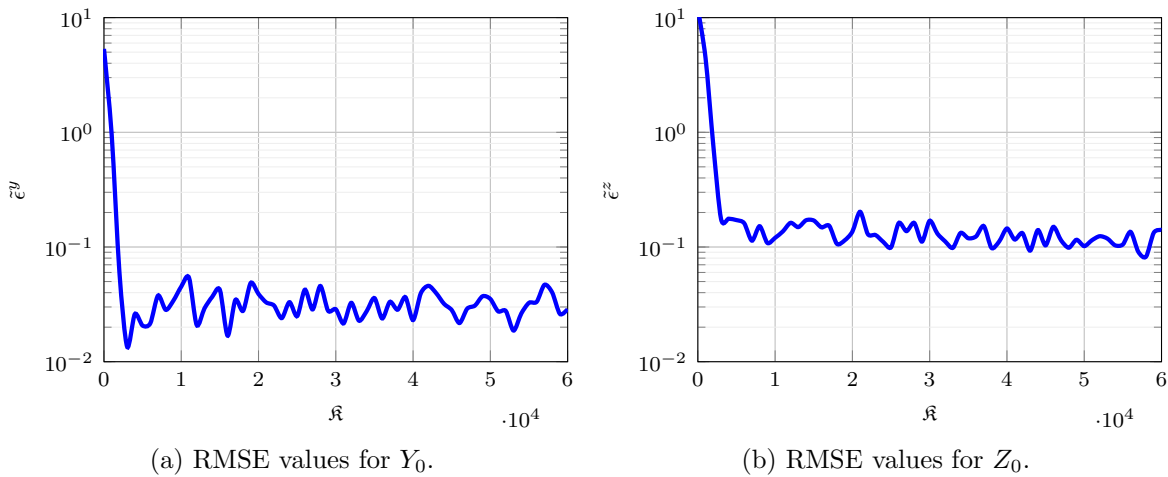
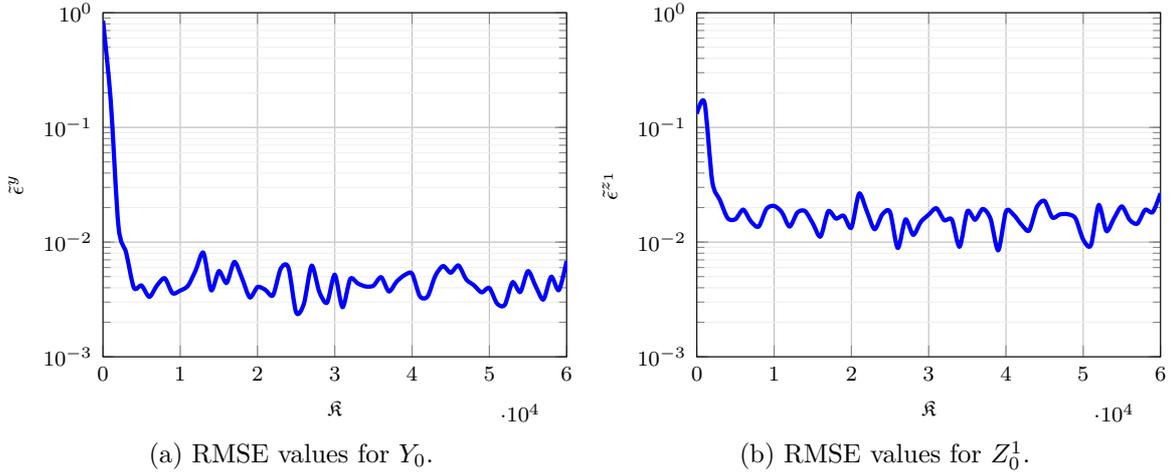
\begin{figure}[htb!]
        \centering
	\begin{subfigure}[h!]{0.48\linewidth}
            \begin{tikzpicture} 
            \begin{axis}[
                xmin = 0, xmax = 60000,
                ymin = 1e-3, ymax = 1,
                xtick distance = 10000,
                ymode=log,
                grid = both,
                width = \textwidth,
                height = 0.8\textwidth,
                major grid style = {lightgray},
                minor grid style = {lightgray!25},
                xlabel = {$\Kf$},
                ylabel = {$\tilde{\epsilon}^y$},]
            \addplot[
                smooth,
                ultra thick,
                blue,
                solid
            ] file {"Figures/Example2/Base/RMSE_Y_0.dat"};
            \end{axis}
            \end{tikzpicture}
		\caption{RMSE values for $Y_0$.}
		\label{fig3a}
	\end{subfigure}
	\begin{subfigure}[h!]{0.48\linewidth}
            \begin{tikzpicture} 
            \begin{axis}[
                xmin = 0, xmax = 60000,
                ymin = 1e-3, ymax = 1,
                xtick distance = 10000,
                ymode=log,
                grid = both,
                width = \textwidth,
                height = 0.8\textwidth,
                major grid style = {lightgray},
                minor grid style = {lightgray!25},
                xlabel = {$\Kf$},
                ylabel = {$\tilde{\epsilon}^{z_1}$},]
            \addplot[
                smooth,
                ultra thick,
                blue,
                solid
            ] file {"Figures/Example2/Base/RMSE_Z_0_1.dat"};
            \end{axis}
            \end{tikzpicture}
		\caption{RMSE values for $Z_0^1$.}
		\label{fig3b}
	\end{subfigure}
	\caption{RMSE values are plotted for Example~\ref{ex2} while increasing $\Kf,$ where $T = 0.25$ and $b = 25$.}
	\label{fig3}
\end{figure}
As $\Kf$ increases, the sum of estimation and optimization errors decreases for both examples as expected. This is because the DBSDE algorithm uses a new sample of size $128$ after each optimization step, and the optimizer tends to perform better with more training data and optimization steps. This reduction in RMSE is evident until around $5000$ optimization steps. However, for $\Kf > 5000$, the RMSE plateaus due to other error sources, such as model and discretization errors, which are higher than the optimization error. Note that the RMSE values for Example~\ref{ex2} are lower than those for Example~\ref{ex1} because the exact solution has a smaller value. To further reduce the optimization error, we use a piecewise constant learning rate approach (PC-LR) with $\alpha = \{\num{1e-2}, \num{3e-3}, \num{1e-3}, \num{3e-4}, \num{1e-4}\}$ and $\Kf = \{ 20000, 30000, 40000, 50000, 60000\}$. We compare the RMSE values using C-LR and PC-LR in Figure~\ref{fig4} for Example~\ref{ex1}. A similar behavior is observed for Example~\ref{ex2}, see Figure~\ref{figx1} in Appendix~\ref{AppendixA}.
\begin{figure}[htb!]
        \centering
	\begin{subfigure}[h!]{0.48\linewidth}
            \pgfplotstableread{"Figures/Example1/Optimization_error/RMSE_Y_0.dat"}{\table}
            \begin{tikzpicture} 
            \begin{axis}[
                xmin = 0, xmax = 60000,
                ymin = 1e-3, ymax = 10,
                xtick distance = 10000,
                ymode=log,
                grid = both,
                width = \textwidth,
                height = 0.8\textwidth,
                major grid style = {lightgray},
                minor grid style = {lightgray!25},
                xlabel = {$\Kf$},
                ylabel = {$\tilde{\epsilon}^y$},
                legend cell align = {left},
                legend pos = north east,
                legend style={nodes={scale=0.7, transform shape}}]
                ]
                \addplot[smooth, ultra thick, blue, solid] table [x = 0, y =1] {\table}; 
                \addplot[smooth, ultra thick, orange, dotted] table [x =0, y = 2] {\table};
                \legend{
                    C-LR, 
                    PC-LR
                }
            \end{axis}
            \end{tikzpicture}
		\caption{RMSE values for $Y_0$.}
		\label{fig4a}
	\end{subfigure}
	\begin{subfigure}[h!]{0.48\linewidth}
        \pgfplotstableread{"Figures/Example1/Optimization_error/RMSE_Z_0.dat"}{\table}
        \begin{tikzpicture} 
            \begin{axis}[
                xmin = 0, xmax = 60000,
                ymin = 1e-3, ymax = 10,
                xtick distance = 10000,
                ymode=log,
                grid = both,
                width = \textwidth,
                height = 0.8\textwidth,
                major grid style = {lightgray},
                minor grid style = {lightgray!25},
                xlabel = {$\Kf$},
                ylabel = {$\tilde{\epsilon}^z$},
                legend cell align = {left},
                legend pos = north east,
                legend style={nodes={scale=0.7, transform shape}}]
                \addplot[smooth, ultra thick, blue, solid] table [x = 0, y =1] {\table}; 
                \addplot[smooth, ultra thick, orange, dotted] table [x =0, y = 2] {\table};
                \legend{
                    C-LR, 
                    PC-LR
                } 
            \end{axis}
            \end{tikzpicture}
		\caption{RMSE values for $Z_0$.}
		\label{fig4b}
	\end{subfigure}
	\caption{RMSE values are plotted for Example~\ref{ex1} using different learning rate approaches, where $T = 1, K = 100, S_0 = 100, a = 0.05, b = 0.2, R = 0.03$ and $\delta = 0$.}
	\label{fig4}
\end{figure}
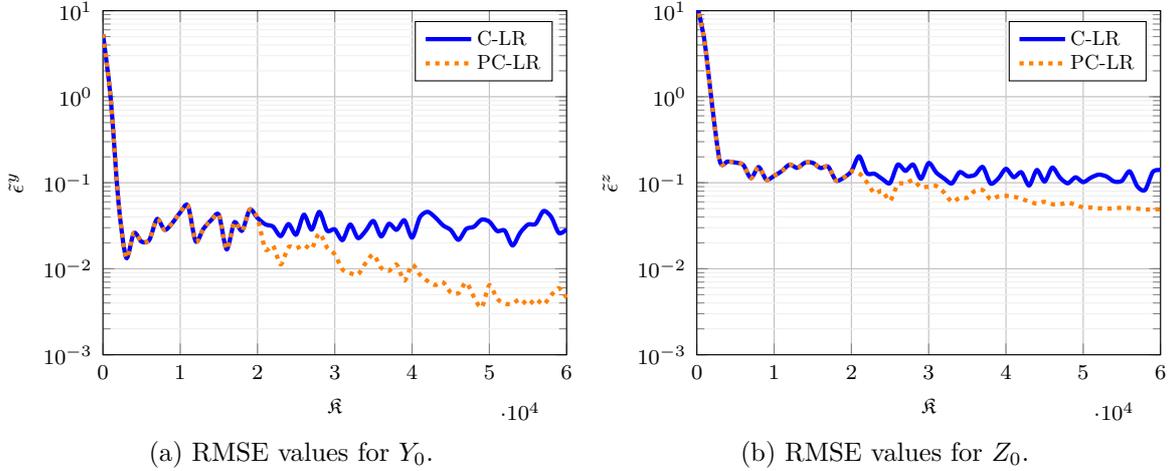

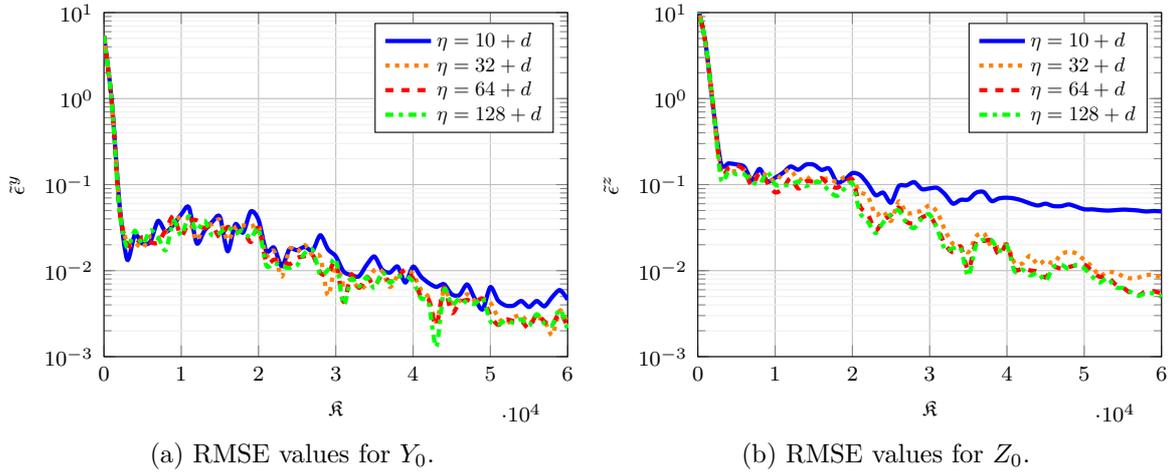
\begin{figure}[htb!]
        \centering
	\begin{subfigure}[h!]{0.48\linewidth}
            \pgfplotstableread{"Figures/Example1/Model_error/RMSE_Y_0.dat"}{\table}
            \begin{tikzpicture} 
            \begin{axis}[
                xmin = 0, xmax = 60000,
                ymin = 1e-3, ymax = 10,
                xtick distance = 10000,
                ymode=log,
                grid = both,
                width = \textwidth,
                height = 0.8\textwidth,
                major grid style = {lightgray},
                minor grid style = {lightgray!25},
                xlabel = {$\Kf$},
                ylabel = {$\tilde{\epsilon}^y$},
                legend cell align = {left},
                legend pos = north east,
                legend style={nodes={scale=0.7, transform shape}}]
                ]
                \addplot[smooth, ultra thick, blue, solid] table [x = 0, y =1] {\table}; 
                \addplot[smooth, ultra thick, orange, dotted] table [x =0, y = 2] {\table};
                \addplot[smooth, ultra thick, red, dashed] table [x =0, y = 3] {\table};
                \addplot[smooth, ultra thick, green, dashdotted] table [x =0, y = 4] {\table};
                
                \legend{
                    $\eta = 10+d$, 
                    $\eta = 32+d$, 
                    $\eta = 64+d$, 
                    $\eta = 128+d$
                } 
            \end{axis}
            \end{tikzpicture}
		\caption{RMSE values for $Y_0$.}
		\label{fig5a}
	\end{subfigure}
	\begin{subfigure}[h!]{0.48\linewidth}
        \pgfplotstableread{"Figures/Example1/Model_error/RMSE_Z_0.dat"}{\table}
        \begin{tikzpicture} 
            \begin{axis}[
                xmin = 0, xmax = 60000,
                ymin = 1e-3, ymax = 10,
                xtick distance = 10000,
                ymode=log,
                grid = both,
                width = \textwidth,
                height = 0.8\textwidth,
                major grid style = {lightgray},
                minor grid style = {lightgray!25},
                xlabel = {$\Kf$},
                ylabel = {$\tilde{\epsilon}^z$},
                legend cell align = {left},
                legend pos = north east,
                legend style={nodes={scale=0.7, transform shape}}]
                \addplot[smooth, ultra thick, blue, solid] table [x = 0, y =1] {\table}; 
                \addplot[smooth, ultra thick, orange, dotted] table [x =0, y = 2] {\table};
                \addplot[smooth, ultra thick, red, dashed] table [x =0, y = 3] {\table};
                \addplot[smooth, ultra thick, green, dashdotted] table [x =0, y = 4] {\table};
                
                \legend{
                    $\eta = 10+d$, 
                    $\eta = 32+d$, 
                    $\eta = 64+d$, 
                    $\eta = 128+d$
                } 
            \end{axis}
            \end{tikzpicture}
		\caption{RMSE values for $Z_0$.}
		\label{fig5b}
	\end{subfigure}
	\caption{RMSE values are plotted for Example~\ref{ex1} using $\eta \in \{ 10+d, 32+d, 64 +d , 128+d \}$, where $T = 1, K = 100, S_0 = 100, a = 0.05, b = 0.2, R = 0.03$ and $\delta = 0$.}
	\label{fig5}
\end{figure}

Next, we consider the model error. To try reduce the model error, one can increase the number of hidden neurons $\eta$ or the number of hidden layers $L$. We report the RMSE values for $\eta \in \{ 10+d, 32+d, 64 +d , 128+d \}$ in Figures~\ref{fig5} and~\ref{fig6} using PC-LR, for Examples~\ref{ex1} and~\ref{ex2}, respectively.
\begin{figure}[htb!]
        \centering
	\begin{subfigure}[h!]{0.48\linewidth}
            \pgfplotstableread{"Figures/Example2/Model_error/RMSE_Y_0.dat"}{\table}
            \begin{tikzpicture} 
            \begin{axis}[
                xmin = 0, xmax = 60000,
                ymin = 1e-4, ymax = 1,
                xtick distance = 10000,
                ymode=log,
                grid = both,
                width = \textwidth,
                height = 0.8\textwidth,
                major grid style = {lightgray},
                minor grid style = {lightgray!25},
                xlabel = {$\Kf$},
                ylabel = {$\tilde{\epsilon}^y$},
                legend cell align = {left},
                legend pos = north east,
                legend style={nodes={scale=0.7, transform shape}}]
                ]
                \addplot[smooth, ultra thick, blue, solid] table [x = 0, y =1] {\table}; 
                \addplot[smooth, ultra thick, orange, dotted] table [x =0, y = 2] {\table};
                \addplot[smooth, ultra thick, red, dashed] table [x =0, y = 3] {\table};
                \addplot[smooth, ultra thick, green, dashdotted] table [x =0, y = 4] {\table};
                
                \legend{
                    $\eta = 10+d$, 
                    $\eta = 32+d$, 
                    $\eta = 64+d$, 
                    $\eta = 128+d$
                } 
            \end{axis}
            \end{tikzpicture}
		\caption{RMSE values for $Y_0$.}
		\label{fig6a}
	\end{subfigure}
	\begin{subfigure}[h!]{0.48\linewidth}
        \pgfplotstableread{"Figures/Example2/Model_error/RMSE_Z_0_1.dat"}{\table}
        \begin{tikzpicture} 
            \begin{axis}[
                xmin = 0, xmax = 60000,
                ymin = 1e-4, ymax = 1,
                xtick distance = 10000,
                ymode=log,
                grid = both,
                width = \textwidth,
                height = 0.8\textwidth,
                major grid style = {lightgray},
                minor grid style = {lightgray!25},
                xlabel = {$\Kf$},
                ylabel = {$\tilde{\epsilon}^{z_1}$},
                legend cell align = {left},
                legend pos = north east,
                legend style={nodes={scale=0.7, transform shape}}]
                \addplot[smooth, ultra thick, blue, solid] table [x = 0, y =1] {\table}; 
                \addplot[smooth, ultra thick, orange, dotted] table [x =0, y = 2] {\table};
                \addplot[smooth, ultra thick, red, dashed] table [x =0, y = 3] {\table};
                \addplot[smooth, ultra thick, green, dashdotted] table [x =0, y = 4] {\table};
                
                \legend{
                    $\eta = 10+d$, 
                    $\eta = 32+d$, 
                    $\eta = 64+d$, 
                    $\eta = 128+d$
                } 
            \end{axis}
            \end{tikzpicture}
		\caption{RMSE values for $Z_0^1$.}
		\label{fig6b}
	\end{subfigure}
	\caption{RMSE values are plotted for Example~\ref{ex2} using $\eta \in \{ 10+d, 32+d, 64 +d , 128+d \}$, where $T = 0.25$ and $b = 25$.}
	\label{fig6}
\end{figure}
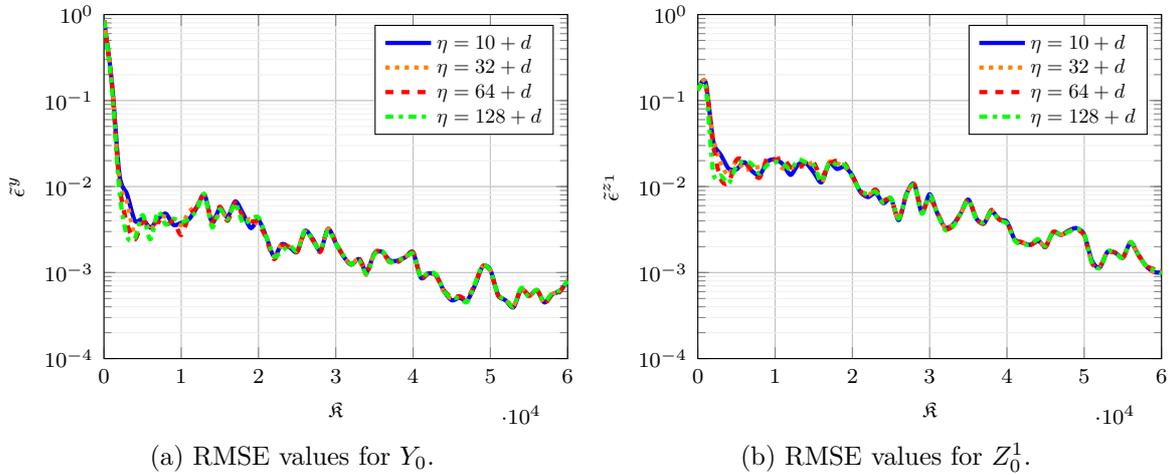We observe that for Example~\ref{ex1} the RMSE decreases when increasing $\eta$, but this trend only persists until $\eta = 64+d$. This is not the case for Example~\ref{ex2}. However, note that $d=50$ in Example~\ref{ex2}, i.e., the starting value of $\eta$ is $60$, which may explain why the RMSE did not decrease as observed in Example~\ref{ex1}.

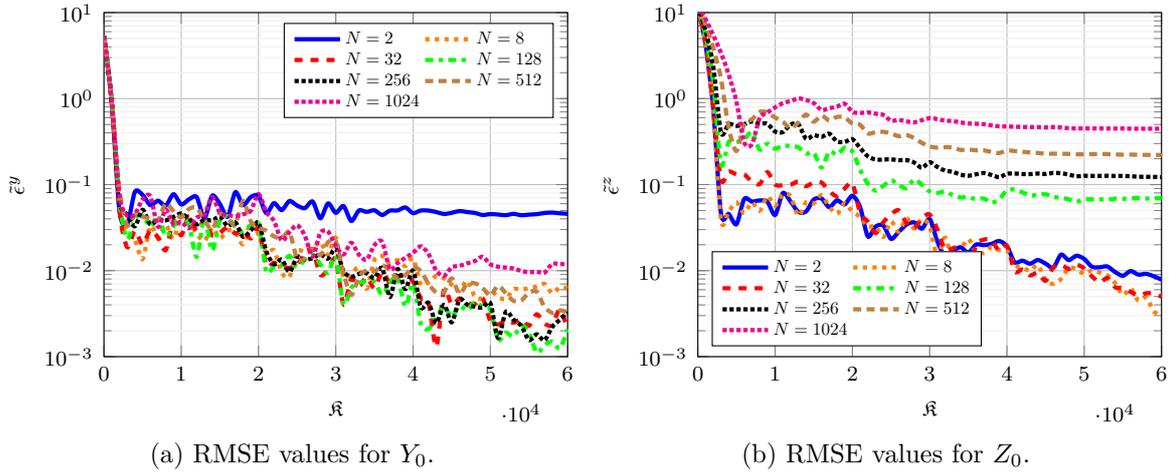
\begin{figure}[htb!]
        \centering
	\begin{subfigure}[h!]{0.48\linewidth}
            \pgfplotstableread{"Figures/Example1/Discretization_error/RMSE_Y_0.dat"}{\table}
            \begin{tikzpicture} 
            \begin{axis}[
                xmin = 0, xmax = 60000,
                ymin = 1e-3, ymax = 10,
                xtick distance = 10000,
                ymode=log,
                grid = both,
                width = \textwidth,
                height = 0.8\textwidth,
                major grid style = {lightgray},
                minor grid style = {lightgray!25},
                xlabel = {$\Kf$},
                ylabel = {$\tilde{\epsilon}^y$},
                legend cell align = {left},
                legend pos = north east,
                legend style={nodes={scale=0.6, transform shape},legend columns=2}]
                ]
                \addplot[smooth, ultra thick, blue, solid] table [x = 0, y =1] {\table}; 
                \addplot[smooth, ultra thick, orange, dotted] table [x =0, y = 2] {\table};
                \addplot[smooth, ultra thick, red, dashed] table [x =0, y = 3] {\table};
                \addplot[smooth, ultra thick, green, dashdotted] table [x =0, y = 4] {\table};
                \addplot[smooth, ultra thick, black, densely dotted] table [x =0, y = 5] {\table};
                \addplot[smooth, ultra thick, brown, densely dashed] table [x =0, y = 6] {\table};
                \addplot[smooth, ultra thick, magenta, densely dotted] table [x =0, y = 7] {\table};
                
                \legend{
                    $N = 2$, 
                    $N = 8$, 
                    $N = 32$, 
                    $N = 128$,
                    $N = 256$,
                    $N = 512$,
                    $N = 1024$
                } 
            \end{axis}
            \end{tikzpicture}
		\caption{RMSE values for $Y_0$.}
		\label{fig7a}
	\end{subfigure}
	\begin{subfigure}[h!]{0.48\linewidth}
        \pgfplotstableread{"Figures/Example1/Discretization_error/RMSE_Z_0.dat"}{\table}
        \begin{tikzpicture} 
            \begin{axis}[
                xmin = 0, xmax = 60000,
                ymin = 1e-3, ymax = 10,
                xtick distance = 10000,
                ymode=log,
                grid = both,
                width = \textwidth,
                height = 0.8\textwidth,
                major grid style = {lightgray},
                minor grid style = {lightgray!25},
                xlabel = {$\Kf$},
                ylabel = {$\tilde{\epsilon}^z$},
                legend cell align = {left},
                legend pos = south west,
                legend style={nodes={scale=0.6, transform shape},legend columns=2}]
                \addplot[smooth, ultra thick, blue, solid] table [x = 0, y =1] {\table}; 
                \addplot[smooth, ultra thick, orange, dotted] table [x =0, y = 2] {\table};
                \addplot[smooth, ultra thick, red, dashed] table [x =0, y = 3] {\table};
                \addplot[smooth, ultra thick, green, dashdotted] table [x =0, y = 4] {\table};
                \addplot[smooth, ultra thick, black, densely dotted] table [x =0, y = 5] {\table};
                \addplot[smooth, ultra thick, brown, densely dashed] table [x =0, y = 6] {\table};
                \addplot[smooth, ultra thick, magenta, densely dotted] table [x =0, y = 7] {\table};                
                \legend{
                    $N = 2$, 
                    $N = 8$, 
                    $N = 32$, 
                    $N = 128$,
                    $N = 256$,
                    $N = 512$,
                    $N = 1024$                    
                } 
            \end{axis}
            \end{tikzpicture}
		\caption{RMSE values for $Z_0$.}
		\label{fig7b}
	\end{subfigure}
	\caption{RMSE values are plotted for Example~\ref{ex1} using $N \in \{ 2, 8, 32 , 128, 256, 512, 1024\}$, where $T = 1, K = 100, S_0 = 100, a = 0.05, b = 0.2, R = 0.03$ and $\delta = 0$.}
	\label{fig7}
\end{figure}

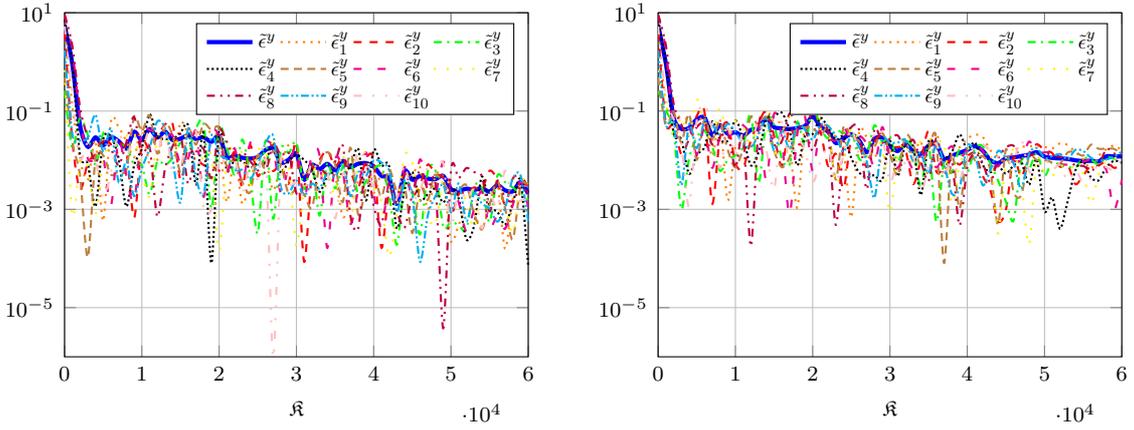
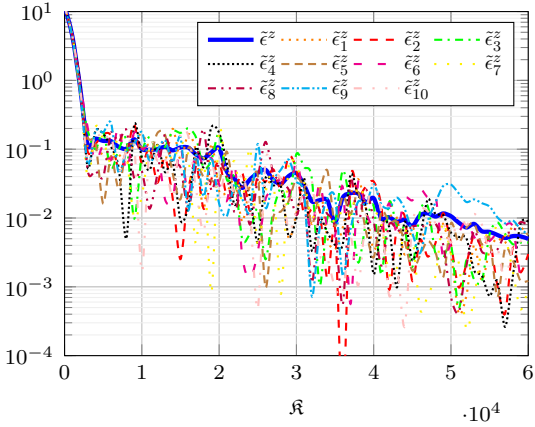
\begin{figure}[htb!]
        \centering
	\begin{subfigure}[h!]{0.48\linewidth}
            \pgfplotstableread{"Figures/Example1/Discretization_error/RMSE_abs_err_Y_0_N_32.dat"}{\table}
            \begin{tikzpicture} 
            \begin{axis}[
                xmin = 0, xmax = 60000,
                ymin = 1e-6, ymax = 10,
                xtick distance = 10000,
                ymode=log,
                grid = both,
                width = \textwidth,
                height = 0.8\textwidth,
                major grid style = {lightgray},
                minor grid style = {lightgray!25},
                xlabel = {$\Kf$},
                legend cell align = {left},
                legend pos = north east,
                legend style={nodes={scale=0.7, transform shape},legend columns=4}]
                ]
                \addplot[smooth, ultra thick, blue, solid] table [x = 0, y =1] {\table}; 
                \addplot[smooth, thick, orange, dotted] table [x =0, y = 2] {\table};
                \addplot[smooth, thick, red, dashed] table [x =0, y = 3] {\table};
                \addplot[smooth, thick, green, dashdotted] table [x =0, y = 4] {\table};
                \addplot[smooth, thick, black, densely dotted] table [x =0, y = 5] {\table};
                \addplot[smooth, thick, brown, densely dashed] table [x =0, y = 6] {\table};
                \addplot[smooth, thick, magenta, loosely dashed] table [x =0, y = 7] {\table};
                \addplot[smooth, thick, yellow, loosely dotted] table [x =0, y = 8] {\table};
                \addplot[smooth, thick, purple, dashdotdotted] table [x =0, y = 9] {\table};
                \addplot[smooth, thick, cyan, densely dashdotdotted] table [x =0, y = 10] {\table};                
                \addplot[smooth, thick, pink, loosely dashdotdotted] table [x =0, y = 11] {\table};                
                \legend{
                    $\tilde{\epsilon}^y$,
                    $\tilde{\epsilon}^y_1$,
                    $\tilde{\epsilon}^y_2$,
                    $\tilde{\epsilon}^y_3$,
                    $\tilde{\epsilon}^y_4$,
                    $\tilde{\epsilon}^y_5$,
                    $\tilde{\epsilon}^y_6$,
                    $\tilde{\epsilon}^y_7$,
                    $\tilde{\epsilon}^y_8$,
                    $\tilde{\epsilon}^y_9$,
                    $\tilde{\epsilon}^y_{10}$
                } 
            \end{axis}
            \end{tikzpicture}
		\caption{RMSE values and the absolute error from each of $Q=10$ DBSDE runs for $Y_0$ with $N = 32$.}
		\label{fig8a}
	\end{subfigure}
	\begin{subfigure}[h!]{0.48\linewidth}
            \pgfplotstableread{"Figures/Example1/Discretization_error/RMSE_abs_err_Y_0_N_1024.dat"}{\table}
            \begin{tikzpicture} 
            \begin{axis}[
                xmin = 0, xmax = 60000,
                ymin = 1e-6, ymax = 10,
                xtick distance = 10000,
                ymode=log,
                grid = both,
                width = \textwidth,
                height = 0.8\textwidth,
                major grid style = {lightgray},
                minor grid style = {lightgray!25},
                xlabel = {$\Kf$},
                legend cell align = {left},
                legend pos = north east,
                legend style={nodes={scale=0.7, transform shape},legend columns=4}]
                ]
                \addplot[smooth, ultra thick, blue, solid] table [x = 0, y =1] {\table}; 
                \addplot[smooth, thick, orange, dotted] table [x =0, y = 2] {\table};
                \addplot[smooth, thick, red, dashed] table [x =0, y = 3] {\table};
                \addplot[smooth, thick, green, dashdotted] table [x =0, y = 4] {\table};
                \addplot[smooth, thick, black, densely dotted] table [x =0, y = 5] {\table};
                \addplot[smooth, thick, brown, densely dashed] table [x =0, y = 6] {\table};
                \addplot[smooth, thick, magenta, loosely dashed] table [x =0, y = 7] {\table};
                \addplot[smooth, thick, yellow, loosely dotted] table [x =0, y = 8] {\table};
                \addplot[smooth, thick, purple, dashdotdotted] table [x =0, y = 9] {\table};
                \addplot[smooth, thick, cyan, densely dashdotdotted] table [x =0, y = 10] {\table};                
                \addplot[smooth, thick, pink, loosely dashdotdotted] table [x =0, y = 11] {\table};                
                \legend{
                    $\tilde{\epsilon}^y$,
                    $\tilde{\epsilon}^y_1$,
                    $\tilde{\epsilon}^y_2$,
                    $\tilde{\epsilon}^y_3$,
                    $\tilde{\epsilon}^y_4$,
                    $\tilde{\epsilon}^y_5$,
                    $\tilde{\epsilon}^y_6$,
                    $\tilde{\epsilon}^y_7$,
                    $\tilde{\epsilon}^y_8$,
                    $\tilde{\epsilon}^y_9$,
                    $\tilde{\epsilon}^y_{10}$
                } 
            \end{axis}
            \end{tikzpicture}
		\caption{RMSE values and the absolute error from each of $Q=10$ DBSDE runs for $Y_0$ with $N = 1024$.}
		\label{fig8b}
	\end{subfigure}
	\begin{subfigure}[h!]{0.48\linewidth}
            \pgfplotstableread{"Figures/Example1/Discretization_error/RMSE_abs_err_Z_0_N_32.dat"}{\table}
            \begin{tikzpicture} 
            \begin{axis}[
                xmin = 0, xmax = 60000,
                ymin = 1e-4, ymax = 10,
                xtick distance = 10000,
                ymode=log,
                grid = both,
                width = \textwidth,
                height = 0.8\textwidth,
                major grid style = {lightgray},
                minor grid style = {lightgray!25},
                xlabel = {$\Kf$},
                legend cell align = {left},
                legend pos = north east,
                legend style={nodes={scale=0.7, transform shape},legend columns=4}]
                ]
                \addplot[smooth, ultra thick, blue, solid] table [x = 0, y =1] {\table}; 
                \addplot[smooth, thick, orange, dotted] table [x =0, y = 2] {\table};
                \addplot[smooth, thick, red, dashed] table [x =0, y = 3] {\table};
                \addplot[smooth, thick, green, dashdotted] table [x =0, y = 4] {\table};
                \addplot[smooth, thick, black, densely dotted] table [x =0, y = 5] {\table};
                \addplot[smooth, thick, brown, densely dashed] table [x =0, y = 6] {\table};
                \addplot[smooth, thick, magenta, loosely dashed] table [x =0, y = 7] {\table};
                \addplot[smooth, thick, yellow, loosely dotted] table [x =0, y = 8] {\table};
                \addplot[smooth, thick, purple, dashdotdotted] table [x =0, y = 9] {\table};
                \addplot[smooth, thick, cyan, densely dashdotdotted] table [x =0, y = 10] {\table};                
                \addplot[smooth, thick, pink, loosely dashdotdotted] table [x =0, y = 11] {\table};                
                \legend{
                    $\tilde{\epsilon}^z$,
                    $\tilde{\epsilon}^z_1$,
                    $\tilde{\epsilon}^z_2$,
                    $\tilde{\epsilon}^z_3$,
                    $\tilde{\epsilon}^z_4$,
                    $\tilde{\epsilon}^z_5$,
                    $\tilde{\epsilon}^z_6$,
                    $\tilde{\epsilon}^z_7$,
                    $\tilde{\epsilon}^z_8$,
                    $\tilde{\epsilon}^z_9$,
                    $\tilde{\epsilon}^z_{10}$
                } 

            \end{axis}
            \end{tikzpicture}
		\caption{RMSE values and the absolute error from each of $Q=10$ DBSDE runs for $Z_0$ with $N = 32$.}
		\label{fig8c}
	\end{subfigure}
	\begin{subfigure}[h!]{0.48\linewidth}
            \pgfplotstableread{"Figures/Example1/Discretization_error/RMSE_abs_err_Z_0_N_1024.dat"}{\table}
            \begin{tikzpicture} 
            \begin{axis}[
                xmin = 0, xmax = 60000,
                ymin = 1e-4, ymax = 10,
                xtick distance = 10000,
                ymode=log,
                grid = both,
                width = \textwidth,
                height = 0.8\textwidth,
                major grid style = {lightgray},
                minor grid style = {lightgray!25},
                xlabel = {$\Kf$},
                legend cell align = {left},
                legend pos = south east,
                legend style={nodes={scale=0.7, transform shape},legend columns=4}]
                ]
                \addplot[smooth, ultra thick, blue, solid] table [x = 0, y =1] {\table}; 
                \addplot[smooth, thick, orange, dotted] table [x =0, y = 2] {\table};
                \addplot[smooth, thick, red, dashed] table [x =0, y = 3] {\table};
                \addplot[smooth, thick, green, dashdotted] table [x =0, y = 4] {\table};
                \addplot[smooth, thick, black, densely dotted] table [x =0, y = 5] {\table};
                \addplot[smooth, thick, brown, densely dashed] table [x =0, y = 6] {\table};
                \addplot[smooth, thick, magenta, loosely dashed] table [x =0, y = 7] {\table};
                \addplot[smooth, thick, yellow, loosely dotted] table [x =0, y = 8] {\table};
                \addplot[smooth, thick, purple, dashdotdotted] table [x =0, y = 9] {\table};
                \addplot[smooth, thick, cyan, densely dashdotdotted] table [x =0, y = 10] {\table};                
                \addplot[smooth, thick, pink, loosely dashdotdotted] table [x =0, y = 11] {\table};                
                \legend{
                    $\tilde{\epsilon}^z$,
                    $\tilde{\epsilon}^z_1$,
                    $\tilde{\epsilon}^z_2$,
                    $\tilde{\epsilon}^z_3$,
                    $\tilde{\epsilon}^z_4$,
                    $\tilde{\epsilon}^z_5$,
                    $\tilde{\epsilon}^z_6$,
                    $\tilde{\epsilon}^z_7$,
                    $\tilde{\epsilon}^z_8$,
                    $\tilde{\epsilon}^z_9$,
                    $\tilde{\epsilon}^z_{10}$
                } 
            \end{axis}
            \end{tikzpicture}
		\caption{RMSE values and the absolute error from each of $Q=10$ DBSDE runs for $Z_0$ with $N = 1024$.}
		\label{fig8d}
	\end{subfigure}
        \caption{RMSE values and the absolute errors from each of $Q=10$ DBSDE runs are plotted for Example~\ref{ex1} using $N \in \{32, 1024\}$, where $T = 1, K = 100, S_0 = 100, a = 0.05, b = 0.2, R = 0.03$ and $\delta = 0$.}
	\label{fig8}
\end{figure}

Finally, we consider the discretization error, which appears to be larger than the model error. To visualize this, we use the PC-LR approach while setting $\eta = 128+d$, and plot the RMSE values in Figure~\ref{fig7} for $N \in \{ 2, 8, 32, 128, 256, 512, 1024\}$ in the case of Example~\ref{ex1}.
 When considering the approximation of $Y_0$, the RMSE decreases with increasing $N$ until $N = 32$, after which it starts to increase. This phenomenon is even more pronounced for the approximation of $Z_0$. It is worth noting that approximating $Z$ is generally more challenging than approximating $Y$ for BSDEs. While a higher value of $N$ can reduce the discretization error, it can also lead to a larger error due to the increased number of DNNs and network parameters to be optimized. Additionally, the propagated errors over time in the DBSDE scheme become larger with a higher value of $N$. A similar behavior is observed for Example~\ref{ex2}, see Figure~\ref{figx2} in Appendix~\ref{AppendixA}. To provide further clarity, we display the RMSE values and the absolute errors $\tilde{\epsilon}_q^y = | Y_{0,q}^{\Delta, \hat{\theta}^m} - Y_0|$ and $\tilde{\epsilon}_q^{z} = | Z_{0,q}^{\Delta, \hat{\theta}^m} - Z_0|$ from the $q$-th run for $N \in \{32, 1024\}$ in Figure~\ref{fig8} for Example~\ref{ex1}.
The variation of absolute errors from different runs around the corresponding RMSE values indicates that the increase in RMSE in Figure~\ref{fig7} for $N > 32$ is caused mainly by the propagated errors (same is observed for Example~\ref{ex2}, see Figure~\ref{figx3} in Appendix~\ref{AppendixA}). Note that choosing $\eta > 128+d$ does not reduce the RMSE. Hence, disentangling the sources of uncertainty is practically challenging, making it essential to quantify the uncertainty of the DBSDE scheme for practical applications.

\input{4_Results_3}

\subsection{Practical implications of the UQ model}
\label{subsec44}
In this section, we study what sources of uncertainty can be captured by our UQ model and we demonstrate its applicability to downstream tasks.

\begin{figure}[htb!]
	\centering
	\begin{subfigure}[h!]{0.48\linewidth}
        \pgfplotstableread{"Figures/Example2/D/Y_0/log_RMSE_STD_STDhat_r_test.dat"}{\table}        
        \begin{tikzpicture} 
        \begin{axis}[
            ymin = -4, ymax = 1,
            xtick distance = 1e-3,
            xmajorgrids={true},
            width = \textwidth,
            height = 0.8\textwidth,
            xlabel = {$\Delta t$},
            legend cell align = {left},
            legend pos = north west,
            legend style={nodes={scale=0.7, transform shape}}
            ]
            \addplot[only marks, black, mark = triangle*, mark options={scale=1}] table [x = 0, y =1] {\table};
            \addplot[only marks, red, mark = square, mark options={scale=1}] table [x =0, y = 2] {\table};
            \addplot[only marks, orange, mark = diamond*, mark options={scale=1}] table [x =0, y = 3] {\table};
            \legend{
            $\log\left( \tilde{\epsilon}^{y,r}(\mathbf{x})\right)$,
            $\log\left( \tilde{\sigma}^{y,r}(\mathbf{x}) \right)$,
            $\log\left( \hat{\sigma}^{y, \hat{\omega}^M_1,r}(\mathbf{x}) \right)$
            } 
            \end{axis}
            \end{tikzpicture}
		\caption{Dataset $\mathcal{D}^{y}$.}
		\label{fig19a}
	\end{subfigure} 
	\begin{subfigure}[h!]{0.48\linewidth}
 		\pgfplotstableread{"Figures/Example2/D/Z_0_1/log_RMSE_STD_STDhat_r_test_Z_0_1.dat"}{\table}        
        \begin{tikzpicture} 
        \begin{axis}[
            ymin = -2, ymax = 2,
            xtick distance = 1e-3,
            xmajorgrids={true},
            width = \textwidth,
            height = 0.8\textwidth,
            xlabel = {$\Delta t$},
            legend cell align = {left},
            legend pos = north east,
            legend style={nodes={scale=0.7, transform shape}}
            ]
            \addplot[only marks, black, mark = triangle*, mark options={scale=1}] table [x = 0, y =1] {\table};
            \addplot[only marks, red, mark = square, mark options={scale=1}] table [x =0, y = 2] {\table};
            \addplot[only marks, orange, mark = diamond*, mark options={scale=1}] table [x =0, y = 3] {\table};
            \legend{
            $\log\left( \tilde{\epsilon}^{z_1,r}(\mathbf{x})\right)$,
            $\log\left( \tilde{\sigma}^{z_1,r}(\mathbf{x}) \right)$,
            $\log\left( \hat{\sigma}^{z_1, \hat{\psi}^M_1,r}(\mathbf{x}) \right)$
            } 
            \end{axis}
            \end{tikzpicture}
		\caption{Dataset $\mathcal{D}^{z}$.}
		\label{fig19b}
	\end{subfigure} 
 	 \caption{Relative RMSE, ensemble STD and estimated STD values from the UQ model for increasing value of $\Delta t$ for $\mathcal{D}$ using the testing sample in Example~\ref{ex2}.}
	\label{fig19}
\end{figure}
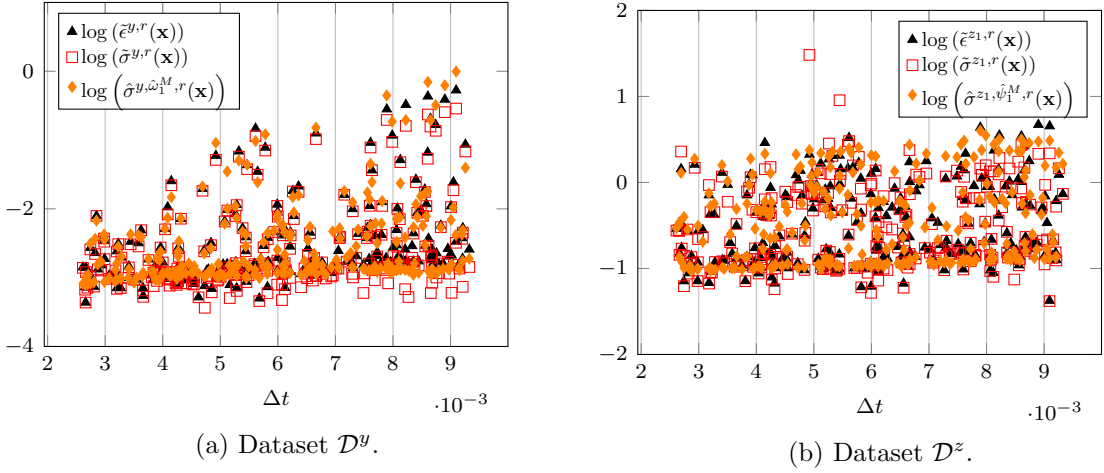

We start by analyzing the sources of uncertainty that our UQ model can effectively capture. It can be expected that the estimated STD captures uncertainty due to the optimization heuristic as well as the uncertainty due to data sampling. However, it is less clear about the uncertainty stemming from the discretization error, as it might bias the approximations provided by the DBSDE scheme. To illustrate the behavior of the relative RMSE, ensemble STD, and estimated STD values across varying $\Delta t$ values, we display these measures in Figure~\ref{fig19} using the testing data in dataset $\mathcal{D}$ from Example~\ref{ex2}.
Note that we use the relative estimated STD from the first trained UQ model (out of our ensemble of 10 models considered for evaluation). As $\Delta t$ decreases, the bias from the discretization error decreases, and the relative estimated STD improves in approximating the relative RMSE. For larger values of $\Delta t$, the bias grows, but the STD also increases. Therefore, the trend of the STD remains consistent with the RMSE, indicating that the relative estimated STD remains reasonable across different values of $\Delta t$ for approximating the relative RMSE. The same is observed for the relative ensemble STD. To measure the strength and direction of the monotonic relationship between the relative RMSE and estimated STD values across $\Delta t$ values, we consider Spearman's rank correlation ($\varsigma$). This metric is calculated as
$$\varsigma\left(\tilde{\epsilon}^{r}(\mathbf{x}),\hat{\sigma}^{r}(\mathbf{x}) \right)=1-\frac{6\sum_{i=1}^{M^{test}} \left( rank(\tilde{\epsilon}^{r}(\mathbf{x}))_i - rank(\hat{\sigma}^{r}(\mathbf{x}))_i \right) }{M^{test}\left(\left(M^{test}\right)^2- 1\right)},$$
where, e.g. $rank(\tilde{\epsilon}^{r}(\mathbf{x}))_i$ is the assigned rank to $\tilde{\epsilon}^{r}(\mathbf{x}_i)$. Note that $\varsigma\left(\tilde{\epsilon}^{r}(\mathbf{x}),\tilde{\sigma}^{r}(\mathbf{x}) \right)$ is calculated similarly. The rank correlation values are displayed in Table~\ref{tab11}.
\begin{table}[h!]
    {\footnotesize
    \begin{center}
      \begin{tabular}{| c | c | c |}\hline
       UQ approach & Rank correlation & Dataset $\mathcal{D}$\\ \hline
      Ensemble for $Y_0$ &  $\varsigma\left( \tilde{\epsilon}^{y,r}(\mathbf{x}), \tilde{\sigma}^{y,r}(\mathbf{x}) \right)$ & $0.9308$\\ \hline
      UQ model for $Y_0$ &  $\varsigma\left(\tilde{\epsilon}^{y,r}(\mathbf{x}),\hat{\sigma}^{y, \hat{\omega}^M_1,r}(\mathbf{x}) \right) $ & $0.8455$\\ \hline
      Ensemble for $Z_0^1$ &  $\varsigma\left( \tilde{\epsilon}^{z_1,r}(\mathbf{x}), \tilde{\sigma}^{z_1,r}(\mathbf{x}) \right) $ & $0.9778$\\ \hline
      UQ model for $Z_0^1$ &  $\varsigma\left(\tilde{\epsilon}^{z_1,r}(\mathbf{x}), \hat{\sigma}^{z_1, \hat{\psi}^M_1,r}(\mathbf{x}) \right)$ & $0.9352$\\ \hline
      \end{tabular}
      \end{center}
    \caption{Rank correlation between the relative RMSE, ensemble STD, and estimated STD values for $\mathcal{D}$ using the testing sample in Example~\ref{ex2}.}
    \label{tab11}  
    }
\end{table}
The high positive rank correlation values indicate that the relative estimated STD from our UQ model can reflect multiple sources of uncertainty, including the uncertainty caused by the discretization error.

Next, we aim to determine whether our UQ model can detect the enhanced performance of the LaDBSDE scheme over the DBSDE scheme for each parameter set, rather than just considering the average performance as shown before. For this purpose, the accuracy score $(acc)$ for the testing sample of datasets $\mathcal{D}$ and $\tilde{\mathcal{D}}$ of Example~\ref{ex2} is considered. We define binary labels to calculate it. For the relative RMSE, we consider
\begin{equation*}
    \ell^{\tilde{\epsilon}^r}(\mathbf{x}_i)=\left\{ 
            \begin{array}{ c l }
    1 & \quad \textrm{if} \,\,  \tilde{\epsilon}^{r, LaDBSDE}(\mathbf{x}_i) < \tilde{\epsilon}^{r, DBSDE}(\mathbf{x}_i), \\
    0                 & \quad \textrm{otherwise,}
  \end{array}
    \right.
\end{equation*}
for the parameter set $\mathbf{x}_i$. Similarly, we define binary labels $\ell^{\tilde{\sigma}^r}(\mathbf{x})$ and  $\ell^{\hat{\sigma}^r}(\mathbf{x}_i)$ for the relative ensemble STD and estimated STD, respectively. The accuracy score between the labels of the relative RMSE and estimated STD values represents the number of parameter sets in which the smallest relative RMSE and estimated STD values are achieved from the same scheme, divided by the total number of parameter sets, i.e.
$$acc\left(\ell^{\tilde{\epsilon}^r}(\mathbf{x}), \ell^{\hat{\sigma}^r}(\mathbf{x}) \right) = \frac{1}{M^{test}}\sum_{i=1}^{M^{test}} \mathds{1}_{\ell^{\tilde{\epsilon}^r}(\mathbf{x}_i)=\ell^{\hat{\sigma}^r}(\mathbf{x}_i)}.$$ 
Similarly, we calculate the accuracy score between the labels of relative RMSE and ensemble STD values $acc\left(\ell^{\tilde{\epsilon}^r}(\mathbf{x}), \ell^{\tilde{\sigma}^r}(\mathbf{x}) \right)$ and report them in Table~\ref{tab12}.
\begin{table}[h!]
    {\footnotesize
    \begin{center}
      \begin{tabular}{| c | c | c |}\hline
       UQ approach & Accuracy score & Datasets $\mathcal{D}$ and $\tilde{\mathcal{D}}$\\ \hline
      Ensemble for $Y_0$ & $acc\left( \ell^{\tilde{\epsilon}^{y,r}}(\mathbf{x}), \ell^{\tilde{\sigma}^{y,r}}(\mathbf{x}) \right)$ & $0.8320$\\ \hline
      UQ model for $Y_0$ & $acc\left( \ell^{\tilde{\epsilon}^{y,r}}(\mathbf{x}), \ell^{\hat{\sigma}^{y,\hat{\omega}^M_1, r}}(\mathbf{x}) \right)$ & $0.7891$\\ \hline
      Ensemble for $Z_0^1$ & $acc\left( \ell^{\tilde{\epsilon}^{z_1,r}}(\mathbf{x}), \ell^{\tilde{\sigma}^{z_1,r}}(\mathbf{x}) \right)$ & $1.0000$\\ \hline
      UQ model for $Z_0^1$ & $acc\left( \ell^{\tilde{\epsilon}^{z_1,r}}(\mathbf{x}), \ell^{\hat{\sigma}^{z_1,\hat{\psi}^M_1, r}}(\mathbf{x}) \right)$ & $1.0000$\\ \hline
      \end{tabular}
      \end{center}
    \caption{Accuracy score between the binary labels of the relative RMSE, ensemble STD, and estimated STD values from $\mathcal{D}$ and $\tilde{\mathcal{D}}$ using the testing sample in Example~\ref{ex2}.}
    \label{tab12}  
    }
\end{table}
The accuracy score of $1$ for $Z_0^1$ implies that the relative estimated STD from our UQ model illustrates enhanced performance when comparing the DBSDE and LaDBSDE schemes across the entire testing sample. This observation is valid only for approximately $80\%$ of the testing sample for $Y_0$.

Moreover, as the RMSE increases due to propagated errors with increasing $N$ (from a certain value of $N$ depending on the parameter set values), it is of interest to determine whether the UQ model can identify the value of $N$ at which the algorithm attains the smallest RMSE based on the estimated STD. To investigate this, we generate a dataset $\mathcal{D}^{\mathbf{N}}$ similar to $\mathcal{D}$ in Table~\ref{tab6} with a fixed maturity $T=0.3$, and each sampled parameter set is solved for $\mathbf{N} = \{ 2, 8, 32, 128\}$. The dataset $\mathcal{D}^\mathbf{N}$ consists of $2560$ parameter sets, resulting in a total number of samples $M = 10240$. We choose $M^{train} = 8192$ and $M^{valid}=M^{test} = 1024$. Note that $\nf=2$ since $\mathbf{x}_i = \left(b_i, N\right)$, $N \in \mathbf{N}$. We use the same hyperparameters for the UQ model as for dataset $\mathcal{D}$ and train only one model. To evaluate the accuracy score, we define the binary multi-label
\begin{equation*}
    \Bell^{\tilde{\epsilon}^r}(\mathbf{x}_i)=\left\{ 
            \begin{array}{ c l }
    \{1, 0, 0, 0 \} & \quad \textrm{if} \,\, N^{min, \tilde{\epsilon}^{r}}(\mathbf{x}_i)  = 2, \\
    \{0, 1, 0, 0 \} & \quad \textrm{if} \,\,  N^{min, \tilde{\epsilon}^{r}}(\mathbf{x}_i) = 8, \\
    \{0, 0, 1, 0 \} & \quad \textrm{if} \,\,  N^{min, \tilde{\epsilon}^{r}}(\mathbf{x}_i) = 32, \\
    \{0, 0, 0, 1 \} & \quad \textrm{if} \,\,  N^{min, \tilde{\epsilon}^{r}}(\mathbf{x}_i) = 128,   
  \end{array}
    \right.
\end{equation*}
for the relative RMSE, where $ N^{min, \tilde{\epsilon}^{r}}(\mathbf{x}_i) = \argmin_{N \in \mathbf{N}} \tilde{\epsilon}^{r}(b_i, N)$. The binary multi-labels for the relative ensemble STD $\Bell^{\tilde{\sigma}^r}(\mathbf{x}_i)$ and estimated STD $\Bell^{\hat{\sigma}^r}(\mathbf{x}_i)$ are defined similarly. The accuracy score values between these multi-labels for the testing sample of $\mathcal{D}^{\mathbf{N}}$ are presented in Table~\ref{tab13}.
\begin{table}[h!]
    {\footnotesize
    \begin{center}
      \begin{tabular}{| c | c | c |}\hline
       UQ approach & Accuracy score & Dataset $\mathcal{D}^{\mathbf{N}}$\\ \hline
      Ensemble for $Y_0$ & $acc\left( \Bell^{\tilde{\epsilon}^{y,r}}(\mathbf{x}), \Bell^{\tilde{\sigma}^{y,r}}(\mathbf{x}) \right)$ & $0.6680$\\ \hline
      UQ model for $Y_0$ & $acc\left( \Bell^{\tilde{\epsilon}^{y,r}}(\mathbf{x}), \Bell^{\hat{\sigma}^{y,\hat{\omega}^M_1, r}}(\mathbf{x}) \right)$ & $0.4805$\\ \hline
      Ensemble for $Z_0^1$ & $acc\left( \Bell^{\tilde{\epsilon}^{z_1,r}}(\mathbf{x}), \Bell^{\tilde{\sigma}^{z_1,r}}(\mathbf{x}) \right)$ & $0.7773$\\ \hline
      UQ model for $Z_0^1$ & $acc\left( \Bell^{\tilde{\epsilon}^{z_1,r}}(\mathbf{x}), \Bell^{\hat{\sigma}^{z_1,\hat{\psi}^M_1, r}}(\mathbf{x}) \right)$ & $0.5469$\\ \hline
      \end{tabular}
      \end{center}
    \caption{Accuracy score between the multi-labels of the relative RMSE, ensemble STD, and estimated STD values from $\mathcal{D}^{\mathbf{N}}$ using the testing sample in Example~\ref{ex2}.}
    \label{tab13}  
    }
\end{table}
We observe that the accuracy score values between the multi-labels of the relative RMSE and estimated STD values $acc\left( \Bell^{\tilde{\epsilon}^{r}}(\mathbf{x}), \Bell^{\hat{\sigma}^{r}}(\mathbf{x}) \right)$ are approximately $0.5$. This indicates that the relative estimated STD correctly predicted the value of N with the smallest relative RMSE for around $50\%$ of the parameter sets in the testing sample.

The accuracy score serves as a restrictive metric, requiring each predicted label ($\Bell^{\tilde{\sigma}^{r}}(\mathbf{x})$ or $\Bell^{\hat{\sigma}^{r}}(\mathbf{x})$) to exactly match the true label $(\Bell^{\tilde{\epsilon}^{r}}(\mathbf{x}))$. Hence, it doesn't tolerate partial errors. For instance, if the $N$ value with the smallest relative RMSE coincides with the one having the second smallest relative estimated STD, the prediction is counted as incorrect. This rigid evaluation fails to consider the order of predicted labels. For this purpose, we consider the mean reciprocal rank $(MRR)$ metric, since there is only one relevant label per sample. It measures the effectiveness of a model in ranking a list of predicted labels based on their relevance to the only true label. In our case, the true label is $N^{min, \tilde{\epsilon}^{r}}(\mathbf{x})$. The predicted ones are denoted by $\mathbf{N}^{sort, \tilde{\sigma}^{r}}(\mathbf{x})$ and $\mathbf{N}^{sort, \hat{\sigma}^{r}}(\mathbf{x})$, the ascending sorted $\mathbf{N}$ values for the parameter set $\mathbf{x}$ based on the value of relative ensemble STD and estimated STD, respectively. Hence, $MRR\left( N^{min, \tilde{\epsilon}^{r}}(\mathbf{x}), \mathbf{N}^{sort, \hat{\sigma}^{r}}(\mathbf{x}) \right)$ is given by
$$MRR\left( N^{min, \tilde{\epsilon}^{r}}(\mathbf{x}), \mathbf{N}^{sort, \hat{\sigma}^{r}}(\mathbf{x}) \right) = \frac{1}{256} \sum_{i=1}^{256} \frac{1}{pos\left( N^{min, \tilde{\epsilon}^{r}}(\mathbf{x}_i), \mathbf{N}^{sort, \hat{\sigma}^{r}}(\mathbf{x}_i)\right)},$$
where $pos\left( N^{min, \tilde{\epsilon}^{r}}(\mathbf{x}_i), \mathbf{N}^{sort, \hat{\sigma}^{r}}(\mathbf{x}_i)\right)$ gives the position where the true label $N^{min, \tilde{\epsilon}^{r}}(\mathbf{x}_i)$ is found in the list of predicted labels $\mathbf{N}^{sort, \hat{\sigma}^{r}}(\mathbf{x}_i)$. The mean reciprocal rank values are reported in Table~\ref{tab14}.
\begin{table}[h!]
    {\footnotesize
    \begin{center}
      \begin{tabular}{| c | c | c |}\hline
       UQ approach & Mean reciprocal rank & Dataset $\mathcal{D}^{\mathbf{N}}$\\ \hline
      Ensemble for $Y_0$ & $MRR\left( N^{min, \tilde{\epsilon}^{y,r}}(\mathbf{x}), \mathbf{N}^{sort, \tilde{\sigma}^{y,r}}(\mathbf{x}) \right)$ & $0.8119$\\ \hline
      UQ model for $Y_0$ & $MRR\left( N^{min, \tilde{\epsilon}^{y,r}}(\mathbf{x}), \mathbf{N}^{sort, \hat{\sigma}^{y, \hat{\omega}^M_1, r}}(\mathbf{x}) \right)$ & $0.6849$\\ \hline
      Ensemble for $Z_0^1$ & $MRR\left( N^{min, \tilde{\epsilon}^{z_1,r}}(\mathbf{x}), \mathbf{N}^{sort, \tilde{\sigma}^{z_1,r}}(\mathbf{x}) \right)$ & $0.8812$\\ \hline
      UQ model for $Z_0^1$ & $MRR\left( N^{min, \tilde{\epsilon}^{z_1,r}}(\mathbf{x}), \mathbf{N}^{sort, \hat{\sigma}^{z_1, \hat{\psi}^M_1, r}}(\mathbf{x}) \right)$ & $0.7614$\\ \hline
      \end{tabular}
      \end{center}
    \caption{Mean reciprocal rank between the $N$ value with the smallest relative RMSE and the ascending sorted $\mathbf{N}$ values based on the relative ensemble STD and estimated STD values from $\mathcal{D}^{\mathbf{N}}$ using the testing sample in Example~\ref{ex2}.}
    \label{tab14}  
    }
\end{table}
We observe that, on average, our UQ model can show that the smallest relative RMSE is achieved for the $N$ value of either the first or second smallest relative estimated STD.

\section{Conclusions}
\label{sec5}
In this work, we investigate the sources of uncertainty in the deep learning-based BSDE schemes and develop a UQ model based on heteroscedastic nonlinear regression to estimate the uncertainty. We apply the UQ model to the pioneering scheme developed in~\cite{weinan2017deep} and the one in~\cite{kapllani2020deep}. The STD of the approximate solution captures the uncertainty, which is usually estimated by performing multiple runs of the algorithm with different datasets. This approach is quite computationally expensive, especially in high-dimensional cases. Our UQ model estimates the STD much cheaper, namely using a single run of the algorithm. Under the assumption of normally distributed errors with zero mean and the STD depending on the parameter set of the discretized BSDE, we employ a DNN to learn two functions that estimate the mean and STD of the approximate solution. The DNN is trained using a dataset of i.i.d. samples, consisting of various parameter sets of the discretized BSDE and their corresponding approximated solutions from a single run of the algorithm. The network parameters are optimized by minimizing the negative log-likelihood. The STD is thus estimated much cheaper. Furthermore, the estimated mean can be leveraged to initialize the algorithm, improving the optimization process. Our numerical results demonstrate that the proposed UQ model provides reliable estimates of the mean and STD of the approximate solution for both considered schemes, even in high-dimensional cases. The estimated STD captures various sources of uncertainty, showcasing its capability in quantifying the uncertainty. Moreover, the UQ model illustrates the improved performance of the LaDBSDE scheme compared to the DBSDE scheme based on the corresponding estimated STD values. Finally, it can also identify the hyperparameters that yield a well-performing scheme.

\begin{appendices}
\section{Impact of the sources of uncertainty for the Burgers type BSDE}
\label{AppendixA}
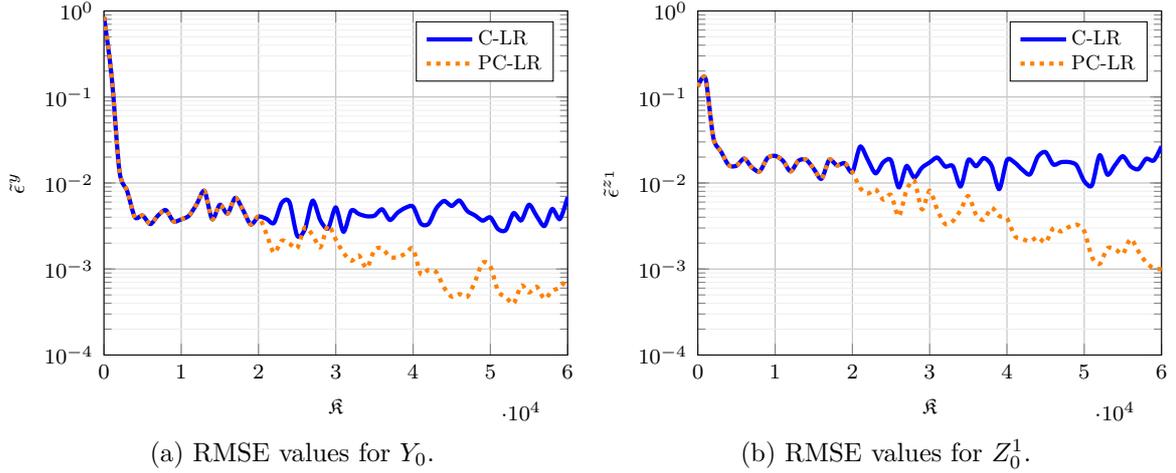
\begin{figure}[htb!]
        \centering
	\begin{subfigure}[h!]{0.48\linewidth}
            \pgfplotstableread{"Figures/Example2/Optimization_error/RMSE_Y_0.dat"}{\table}
            \begin{tikzpicture} 
            \begin{axis}[
                xmin = 0, xmax = 60000,
                ymin = 1e-4, ymax = 1,
                xtick distance = 10000,
                ymode=log,
                grid = both,
                width = \textwidth,
                height = 0.8\textwidth,
                major grid style = {lightgray},
                minor grid style = {lightgray!25},
                xlabel = {$\Kf$},
                ylabel = {$\tilde{\epsilon}^y$},
                legend cell align = {left},
                legend pos = north east,
                legend style={nodes={scale=0.7, transform shape}}]
                ]
                \addplot[smooth, ultra thick, blue, solid] table [x = 0, y =1] {\table}; 
                \addplot[smooth, ultra thick, orange, dotted] table [x =0, y = 2] {\table};
                \legend{
                    C-LR, 
                    PC-LR
                }
            \end{axis}
            \end{tikzpicture}
		\caption{RMSE values for $Y_0$.}
		\label{figx1a}
	\end{subfigure}
	\begin{subfigure}[h!]{0.48\linewidth}
        \pgfplotstableread{"Figures/Example2/Optimization_error/RMSE_Z_0_1.dat"}{\table}
        \begin{tikzpicture} 
            \begin{axis}[
                xmin = 0, xmax = 60000,
                ymin = 1e-4, ymax = 1,
                xtick distance = 10000,
                ymode=log,
                grid = both,
                width = \textwidth,
                height = 0.8\textwidth,
                major grid style = {lightgray},
                minor grid style = {lightgray!25},
                xlabel = {$\Kf$},
                ylabel = {$\tilde{\epsilon}^{z_1}$},
                legend cell align = {left},
                legend pos = north east,
                legend style={nodes={scale=0.7, transform shape}}]
                \addplot[smooth, ultra thick, blue, solid] table [x = 0, y =1] {\table}; 
                \addplot[smooth, ultra thick, orange, dotted] table [x =0, y = 2] {\table};
                \legend{
                    C-LR, 
                    PC-LR
                } 
            \end{axis}
            \end{tikzpicture}
		\caption{RMSE values for $Z_0^1$.}
		\label{figx1b}
	\end{subfigure}
	\caption{RMSE values are plotted for Example~\ref{ex2} using different learning rate approaches, where $T = 0.25$ and $b = 25$.}
	\label{figx1}
\end{figure}

\begin{figure}[htb!]
        \centering
	\begin{subfigure}[h!]{0.48\linewidth}
            \pgfplotstableread{"Figures/Example2/Discretization_error/RMSE_Y_0.dat"}{\table}
            \begin{tikzpicture} 
            \begin{axis}[
                xmin = 0, xmax = 60000,
                ymin = 1e-4, ymax = 1,
                xtick distance = 10000,
                ymode=log,
                grid = both,
                width = \textwidth,
                height = 0.8\textwidth,
                major grid style = {lightgray},
                minor grid style = {lightgray!25},
                xlabel = {$\Kf$},
                ylabel = {$\tilde{\epsilon}^y$},
                legend cell align = {left},
                legend pos = north east,
                legend style={nodes={scale=0.6, transform shape},legend columns=2}]
                ]
                \addplot[smooth, ultra thick, blue, solid] table [x = 0, y =1] {\table}; 
                \addplot[smooth, ultra thick, orange, dotted] table [x =0, y = 2] {\table};
                \addplot[smooth, ultra thick, red, dashed] table [x =0, y = 3] {\table};
                \addplot[smooth, ultra thick, green, dashdotted] table [x =0, y = 4] {\table};
                \addplot[smooth, ultra thick, black, densely dotted] table [x =0, y = 5] {\table};
                \addplot[smooth, ultra thick, brown, densely dashed] table [x =0, y = 6] {\table};
                \addplot[smooth, ultra thick, magenta, densely dotted] table [x =0, y = 7] {\table};
                
                \legend{
                    $N = 2$, 
                    $N = 8$, 
                    $N = 32$, 
                    $N = 128$,
                    $N = 256$,
                    $N = 512$,
                    $N = 1024$
                } 
            \end{axis}
            \end{tikzpicture}
		\caption{RMSE values for $Y_0$.}
		\label{figx2a}
	\end{subfigure}
	\begin{subfigure}[h!]{0.48\linewidth}
        \pgfplotstableread{"Figures/Example2/Discretization_error/RMSE_Z_0_1.dat"}{\table}
        \begin{tikzpicture} 
            \begin{axis}[
                xmin = 0, xmax = 60000,
                ymin = 1e-4, ymax = 1,
                xtick distance = 10000,
                ymode=log,
                grid = both,
                width = \textwidth,
                height = 0.8\textwidth,
                major grid style = {lightgray},
                minor grid style = {lightgray!25},
                xlabel = {$\Kf$},
                ylabel = {$\tilde{\epsilon}^{z_1}$},
                legend cell align = {left},
                legend pos = south west,
                legend style={nodes={scale=0.6, transform shape},legend columns=2}]
                \addplot[smooth, ultra thick, blue, solid] table [x = 0, y =1] {\table}; 
                \addplot[smooth, ultra thick, orange, dotted] table [x =0, y = 2] {\table};
                \addplot[smooth, ultra thick, red, dashed] table [x =0, y = 3] {\table};
                \addplot[smooth, ultra thick, green, dashdotted] table [x =0, y = 4] {\table};
                \addplot[smooth, ultra thick, black, densely dotted] table [x =0, y = 5] {\table};
                \addplot[smooth, ultra thick, brown, densely dashed] table [x =0, y = 6] {\table};
                \addplot[smooth, ultra thick, magenta, densely dotted] table [x =0, y = 7] {\table};                
                \legend{
                    $N = 2$, 
                    $N = 8$, 
                    $N = 32$, 
                    $N = 128$,
                    $N = 256$,
                    $N = 512$,
                    $N = 1024$                    
                } 
            \end{axis}
            \end{tikzpicture}
		\caption{RMSE values for $Z_0^1$.}
		\label{figx2b}
	\end{subfigure}
	\caption{RMSE values are plotted for Example~\ref{ex1} using $N \in \{ 2, 8, 32, 128, 256, 512, 1024\}$, where $T = 0.25$ and $b = 25$.}
	\label{figx2}
\end{figure}
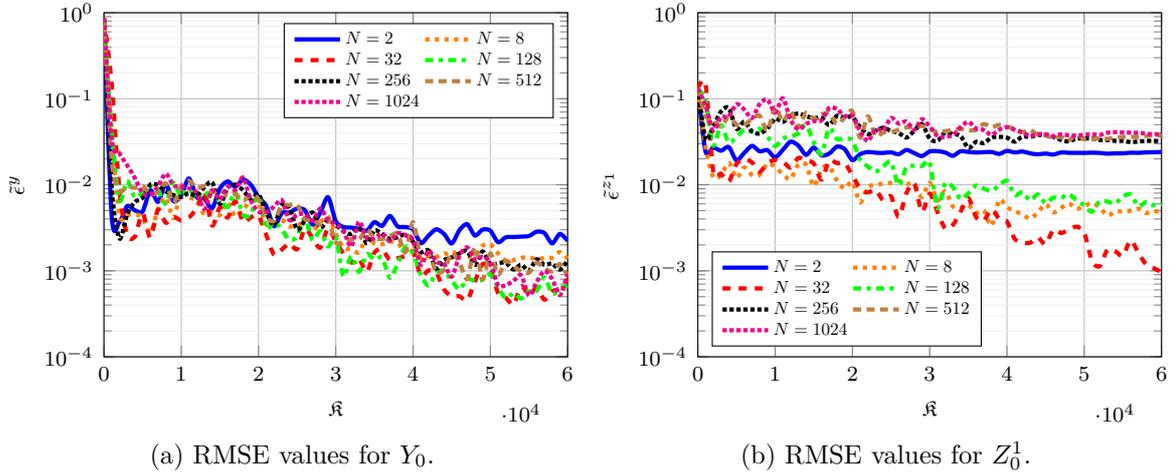

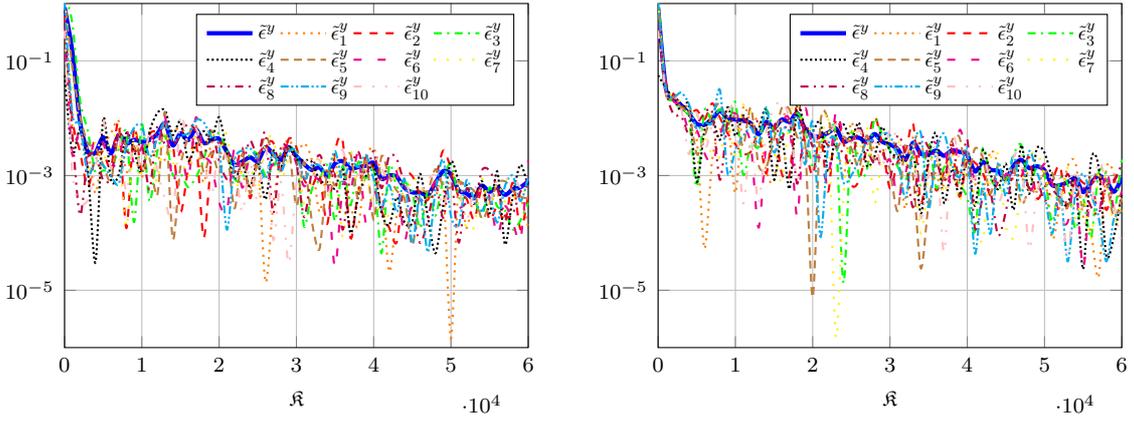
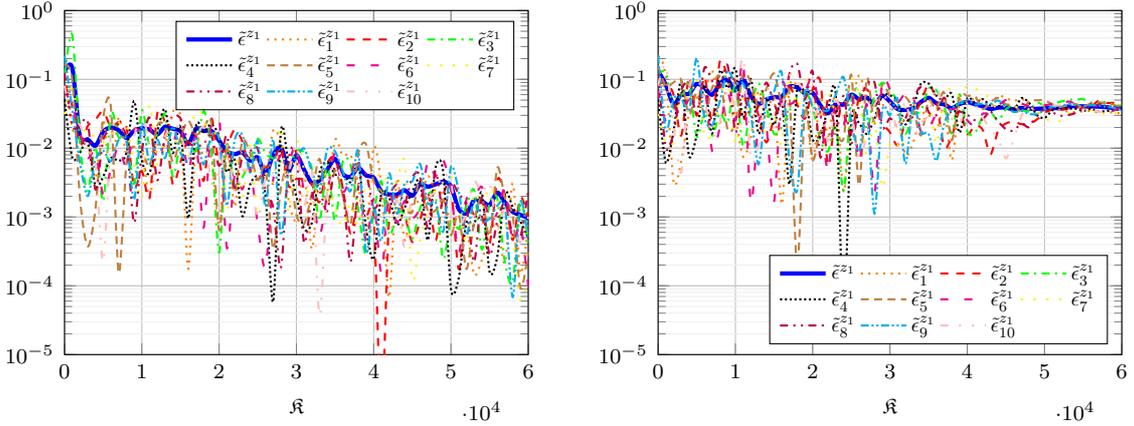
\begin{figure}[htb!]
        \centering
	\begin{subfigure}[h!]{0.48\linewidth}
            \pgfplotstableread{"Figures/Example2/Discretization_error/RMSE_abs_err_Y_0_N_32.dat"}{\table}
            \begin{tikzpicture} 
            \begin{axis}[
                xmin = 0, xmax = 60000,
                ymin = 1e-6, ymax = 1,
                xtick distance = 10000,
                ymode=log,
                grid = both,
                width = \textwidth,
                height = 0.8\textwidth,
                major grid style = {lightgray},
                minor grid style = {lightgray!25},
                xlabel = {$\Kf$},
                legend cell align = {left},
                legend pos = north east,
                legend style={nodes={scale=0.7, transform shape},legend columns=4}]
                ]
                \addplot[smooth, ultra thick, blue, solid] table [x = 0, y =1] {\table}; 
                \addplot[smooth, thick, orange, dotted] table [x =0, y = 2] {\table};
                \addplot[smooth, thick, red, dashed] table [x =0, y = 3] {\table};
                \addplot[smooth, thick, green, dashdotted] table [x =0, y = 4] {\table};
                \addplot[smooth, thick, black, densely dotted] table [x =0, y = 5] {\table};
                \addplot[smooth, thick, brown, densely dashed] table [x =0, y = 6] {\table};
                \addplot[smooth, thick, magenta, loosely dashed] table [x =0, y = 7] {\table};
                \addplot[smooth, thick, yellow, loosely dotted] table [x =0, y = 8] {\table};
                \addplot[smooth, thick, purple, dashdotdotted] table [x =0, y = 9] {\table};
                \addplot[smooth, thick, cyan, densely dashdotdotted] table [x =0, y = 10] {\table};                
                \addplot[smooth, thick, pink, loosely dashdotdotted] table [x =0, y = 11] {\table};                
                \legend{
                    $\tilde{\epsilon}^y$,
                    $\tilde{\epsilon}^y_1$,
                    $\tilde{\epsilon}^y_2$,
                    $\tilde{\epsilon}^y_3$,
                    $\tilde{\epsilon}^y_4$,
                    $\tilde{\epsilon}^y_5$,
                    $\tilde{\epsilon}^y_6$,
                    $\tilde{\epsilon}^y_7$,
                    $\tilde{\epsilon}^y_8$,
                    $\tilde{\epsilon}^y_9$,
                    $\tilde{\epsilon}^y_{10}$
                } 
            \end{axis}
            \end{tikzpicture}
		\caption{RMSE values and the absolute error for each of $Q=10$ DBSDE runs for $Y_0$ with $N = 32$.}
		\label{figx3a}
	\end{subfigure}
	\begin{subfigure}[h!]{0.48\linewidth}
            \pgfplotstableread{"Figures/Example2/Discretization_error/RMSE_abs_err_Y_0_N_1024.dat"}{\table}
            \begin{tikzpicture} 
            \begin{axis}[
                xmin = 0, xmax = 60000,
                ymin = 1e-6, ymax = 1,
                xtick distance = 10000,
                ymode=log,
                grid = both,
                width = \textwidth,
                height = 0.8\textwidth,
                major grid style = {lightgray},
                minor grid style = {lightgray!25},
                xlabel = {$\Kf$},
                legend cell align = {left},
                legend pos = north east,
                legend style={nodes={scale=0.7, transform shape},legend columns=4}]
                ]
                \addplot[smooth, ultra thick, blue, solid] table [x = 0, y =1] {\table}; 
                \addplot[smooth, thick, orange, dotted] table [x =0, y = 2] {\table};
                \addplot[smooth, thick, red, dashed] table [x =0, y = 3] {\table};
                \addplot[smooth, thick, green, dashdotted] table [x =0, y = 4] {\table};
                \addplot[smooth, thick, black, densely dotted] table [x =0, y = 5] {\table};
                \addplot[smooth, thick, brown, densely dashed] table [x =0, y = 6] {\table};
                \addplot[smooth, thick, magenta, loosely dashed] table [x =0, y = 7] {\table};
                \addplot[smooth, thick, yellow, loosely dotted] table [x =0, y = 8] {\table};
                \addplot[smooth, thick, purple, dashdotdotted] table [x =0, y = 9] {\table};
                \addplot[smooth, thick, cyan, densely dashdotdotted] table [x =0, y = 10] {\table};                
                \addplot[smooth, thick, pink, loosely dashdotdotted] table [x =0, y = 11] {\table};                
                \legend{
                    $\tilde{\epsilon}^y$,
                    $\tilde{\epsilon}^y_1$,
                    $\tilde{\epsilon}^y_2$,
                    $\tilde{\epsilon}^y_3$,
                    $\tilde{\epsilon}^y_4$,
                    $\tilde{\epsilon}^y_5$,
                    $\tilde{\epsilon}^y_6$,
                    $\tilde{\epsilon}^y_7$,
                    $\tilde{\epsilon}^y_8$,
                    $\tilde{\epsilon}^y_9$,
                    $\tilde{\epsilon}^y_{10}$
                } 
            \end{axis}
            \end{tikzpicture}
		\caption{RMSE values and the absolute error for each of $Q=10$ DBSDE runs for $Y_0$ with $N = 1024$.}
		\label{figx3b}
	\end{subfigure}
	\begin{subfigure}[h!]{0.48\linewidth}
            \pgfplotstableread{"Figures/Example2/Discretization_error/RMSE_abs_err_Z_0_1_N_32.dat"}{\table}
            \begin{tikzpicture} 
            \begin{axis}[
                xmin = 0, xmax = 60000,
                ymin = 1e-5, ymax = 1,
                xtick distance = 10000,
                ymode=log,
                grid = both,
                width = \textwidth,
                height = 0.8\textwidth,
                major grid style = {lightgray},
                minor grid style = {lightgray!25},
                xlabel = {$\Kf$},
                legend cell align = {left},
                legend pos = north east,
                legend style={nodes={scale=0.7, transform shape},legend columns=4}]
                ]
                \addplot[smooth, ultra thick, blue, solid] table [x = 0, y =1] {\table}; 
                \addplot[smooth, thick, orange, dotted] table [x =0, y = 2] {\table};
                \addplot[smooth, thick, red, dashed] table [x =0, y = 3] {\table};
                \addplot[smooth, thick, green, dashdotted] table [x =0, y = 4] {\table};
                \addplot[smooth, thick, black, densely dotted] table [x =0, y = 5] {\table};
                \addplot[smooth, thick, brown, densely dashed] table [x =0, y = 6] {\table};
                \addplot[smooth, thick, magenta, loosely dashed] table [x =0, y = 7] {\table};
                \addplot[smooth, thick, yellow, loosely dotted] table [x =0, y = 8] {\table};
                \addplot[smooth, thick, purple, dashdotdotted] table [x =0, y = 9] {\table};
                \addplot[smooth, thick, cyan, densely dashdotdotted] table [x =0, y = 10] {\table};                
                \addplot[smooth, thick, pink, loosely dashdotdotted] table [x =0, y = 11] {\table};                
                \legend{
                    $\tilde{\epsilon}^{z_1}$,
                    $\tilde{\epsilon}^{z_1}_1$,
                    $\tilde{\epsilon}^{z_1}_2$,
                    $\tilde{\epsilon}^{z_1}_3$,
                    $\tilde{\epsilon}^{z_1}_4$,
                    $\tilde{\epsilon}^{z_1}_5$,
                    $\tilde{\epsilon}^{z_1}_6$,
                    $\tilde{\epsilon}^{z_1}_7$,
                    $\tilde{\epsilon}^{z_1}_8$,
                    $\tilde{\epsilon}^{z_1}_9$,
                    $\tilde{\epsilon}^{z_1}_{10}$
                } 

            \end{axis}
            \end{tikzpicture}
		\caption{RMSE values and the absolute error for each of $Q=10$ DBSDE runs for $Z_0^1$ with $N = 32$.}
		\label{figx3c}
	\end{subfigure}
	\begin{subfigure}[h!]{0.48\linewidth}
            \pgfplotstableread{"Figures/Example2/Discretization_error/RMSE_abs_err_Z_0_1_N_1024.dat"}{\table}
            \begin{tikzpicture} 
            \begin{axis}[
                xmin = 0, xmax = 60000,
                ymin = 1e-5, ymax = 1,
                xtick distance = 10000,
                ymode=log,
                grid = both,
                width = \textwidth,
                height = 0.8\textwidth,
                major grid style = {lightgray},
                minor grid style = {lightgray!25},
                xlabel = {$\Kf$},
                legend cell align = {left},
                legend pos = south east,
                legend style={nodes={scale=0.7, transform shape},legend columns=4}]
                ]
                \addplot[smooth, ultra thick, blue, solid] table [x = 0, y =1] {\table}; 
                \addplot[smooth, thick, orange, dotted] table [x =0, y = 2] {\table};
                \addplot[smooth, thick, red, dashed] table [x =0, y = 3] {\table};
                \addplot[smooth, thick, green, dashdotted] table [x =0, y = 4] {\table};
                \addplot[smooth, thick, black, densely dotted] table [x =0, y = 5] {\table};
                \addplot[smooth, thick, brown, densely dashed] table [x =0, y = 6] {\table};
                \addplot[smooth, thick, magenta, loosely dashed] table [x =0, y = 7] {\table};
                \addplot[smooth, thick, yellow, loosely dotted] table [x =0, y = 8] {\table};
                \addplot[smooth, thick, purple, dashdotdotted] table [x =0, y = 9] {\table};
                \addplot[smooth, thick, cyan, densely dashdotdotted] table [x =0, y = 10] {\table};                
                \addplot[smooth, thick, pink, loosely dashdotdotted] table [x =0, y = 11] {\table};                
                \legend{
                    $\tilde{\epsilon}^{z_1}$,
                    $\tilde{\epsilon}^{z_1}_1$,
                    $\tilde{\epsilon}^{z_1}_2$,
                    $\tilde{\epsilon}^{z_1}_3$,
                    $\tilde{\epsilon}^{z_1}_4$,
                    $\tilde{\epsilon}^{z_1}_5$,
                    $\tilde{\epsilon}^{z_1}_6$,
                    $\tilde{\epsilon}^{z_1}_7$,
                    $\tilde{\epsilon}^{z_1}_8$,
                    $\tilde{\epsilon}^{z_1}_9$,
                    $\tilde{\epsilon}^{z_1}_{10}$
                } 
            \end{axis}
            \end{tikzpicture}
		\caption{RMSE values and the absolute error for each of $Q=10$ DBSDE runs for $Z_0^1$ with $N = 1024$.}
		\label{figx3d}
	\end{subfigure}
        \caption{RMSE values and the absolute errors from each of $Q=10$ DBSDE runs are plotted for Example~\ref{ex2} using $N \in \{32, 1024\}$, where $T = 0.25$ and $b = 25$.}
	\label{figx3}
\end{figure}

In this section, we visualize the effect of different errors on the RMSE for Example~\ref{ex2}. The impact of the optimization error is shown in Figure~\ref{figx1} using C-LR and PC-LR approaches.
For the discretization error, see Figure~\ref{figx2}.
The effect of the optimization error and propagated errors over time are displayed in Figure~\ref{figx3}. 

\section{Normality assumption of the error distribution}
\label{AppendixB}
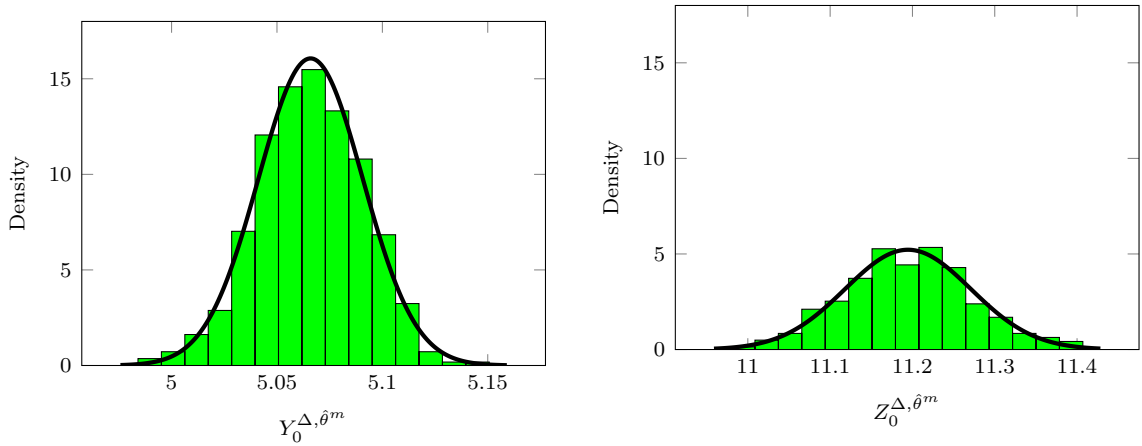
\begin{figure}[htb!]
        \centering
	\begin{subfigure}[h!]{0.48\linewidth}
            \pgfplotstableread{"Figures/Example1/Distribution/BS_Y_0_approx.dat"}{\table}
            \begin{tikzpicture} 
            \begin{axis}[
                ymin = 0, ymax = 18,
                width = \textwidth,
                height = 0.8\textwidth,
                xlabel = {$Y_0^{\Delta, \hat{\theta}^m}$},
                ylabel = {Density}
                ]
                \addplot[ybar,fill=green,
                hist={
                bins=15,
                density
                }
                ] table [y = 0] {\table};
                \addplot[
                smooth,
                ultra thick,
                black,
                solid
                ] file {"Figures/Example1/Distribution/BS_Y_0_fit.dat"};
            \end{axis}
            \end{tikzpicture}
		\caption{Empirical distribution of $Y_0^{\Delta, \hat{\theta}^m}$, where the calculated parameters are $\left(\Tilde{\mu}^y, \Tilde{\sigma}^y\right) =\left(5.0659, 0.0248\right)$.}
		\label{figx4a}
	\end{subfigure}
	\begin{subfigure}[h!]{0.48\linewidth}
        \pgfplotstableread{"Figures/Example1/Distribution/BS_Z_0_approx.dat"}{\table}
            \begin{tikzpicture} 
            \begin{axis}[
                ymin = 0, ymax = 18,
                width = \textwidth,
                height = 0.8\textwidth,
                xlabel = {$Z_0^{\Delta, \hat{\theta}^m}$},
                ylabel = {Density}
                ]
                \addplot[ybar,fill=green,
                hist={
                bins=15,
                density
                }
                ] table [y = 0] {\table};
                \addplot[
                smooth,
                ultra thick,
                black,
                solid
                ] file {"Figures/Example1/Distribution/BS_Z_0_fit.dat"};
            \end{axis}
            \end{tikzpicture}
            
        \caption{Empirical distribution of $Z_0^{\Delta, \hat{\theta}^m}$, where the calculated parameters are $\left(\Tilde{\mu}^z, \Tilde{\sigma}^z\right) =\left(11.1946, 0.0764\right)$.}
		\label{figx4b}
	\end{subfigure}
	\caption{Empirical distribution of the approximate solution~\eqref{eq13} in Example~\ref{ex1} for parameter set $T = 0.33, K = 100, S_0 = 100, a = 0.05, b = 0.2, R = 0.03$ and $\delta = 0$. The curve represents a fitted normal distribution to the data.}
	\label{figx4}
\end{figure}

In this section, we conduct a test to assess the normality of the error distribution in~\eqref{eq13} for Example~\ref{ex1}. For the parameter values $T = 0.33, K = 100, S_0 = 100, a = 0.05, b = 0.2, R = 0.03$ and $\delta = 0$, the exact solution is $(Y_0, Z_0) = (5.0679, 11.1420)$. Using $N = 16$, $\Kf=30000$, $\alpha = \num{1e-2}$ and conducting $Q=500$ independent runs of the DBSDE algorithm, we display the empirical distribution of the approximations in Figure~\ref{figx4}.
The observed empirical distributions exhibit a Gaussian shape. To further assess the normality, we perform the Shapiro-Wilk~\cite{shapiro1965analysis} and D'Agostino and Pearson's~\cite{d1973tests} tests for the assessment of normality. The $p$-values obtained from these tests are presented in Table~\ref{tabx1}.
\begin{table}[h!]
    {\footnotesize
    \begin{center}
      \begin{tabular}{| c | c | c |}
        \cline{2-3}
      \multicolumn{1}{c|}{} & Shapiro-Wilk & D'Agostino-Pearson\\ \hline
         $Y_0^{\Delta, \hat{\theta}^m}$ & $0.5014$ & $0.5368$\\ \hline
         $Z_0^{\Delta, \hat{\theta}^m}$ & $0.4391$ & $0.6489$\\ \hline
        \end{tabular}
      \end{center}
    \caption{$p$-value of the statistical tests in Example~\ref{ex1} for parameter set $T = 0.33, K = 100, S_0 = 100, a = 0.05, b = 0.2, R = 0.03$ and $\delta = 0$.}
    \label{tabx1}  
    }
\end{table}
With a significance level of $0.05$, we conclude that there is no evidence to reject the assumption of normal distribution for the approximate solutions $Y_0^{\Delta,\hat{\theta}^m}$ and $Z_0^{\Delta,\hat{\theta}^m}$.
\end{appendices}

\bibliographystyle{siamplain}
\bibliography{bibfile}

\end{document}